\documentclass[11pt]{amsart}

\usepackage{amsmath}
\usepackage{amssymb}
\usepackage{graphicx}
\newtheorem{theorem}{Theorem}[section]
\newtheorem{corollary}[theorem]{Corollary}
\newtheorem{lemma}[theorem]{Lemma}
\newtheorem{proposition}[theorem]{Proposition}
\theoremstyle{definition}
\newtheorem{definition}[theorem]{Definition}

\newcommand{\argmin}{\operatornamewithlimits{argmin}} %argmax

\newcommand{\RR}{\mathbb{R}}
\newcommand{\Rn}{\mathbb{R}^n}
\newcommand{\Rm}{\mathbb{R}^m}

\newcommand{\be}{\begin{enumerate}}
\newcommand{\ee}{\end{enumerate}}
\newcommand{\bi}{\begin{itemize}}
\newcommand{\ei}{\end{itemize}}
\newcommand{\bc}{\begin{center}}
\newcommand{\ec}{\end{center}}
\newcommand{\beq}{\begin{equation}}
\newcommand{\eeq}{\end{equation}}

\newcommand{\bth}{\begin{theorem}}
\newcommand{\bcor}{\begin{corollary}}
\newcommand{\ecor}{\end{corollary}}
\newcommand{\bax}{\begin{axiom}}
\newcommand{\eax}{\end{axiom}}
\newcommand{\bdf}{\begin{definition}}
\newcommand{\edf}{\end{definition}}
\newcommand{\blem}{\begin{lemma}}
\newcommand{\elem}{\end{lemma}}
\newtheorem{exa}{Example}
\newcommand{\bex}{\begin{exa}}
\newcommand{\eex}{\end{exa}}

\def\tinf{\rightarrow\infty}
\def\tzer{\rightarrow 0}

\def\diag{\mbox{diag}}

\def\lb{\left (}
\def\rb{\right )}
\def\rsb{\right ]}
\def\lsb{\left [}

\def\eps{\varepsilon}
\def\els{\varepsilon_{ls}}

\def\rank{\mbox{rank }}

\def\eee {\eps_0,\eps_1,\els}
\newfont{\gothic}{eufm10}

\def\l{\lambda}
\def\lm{\lambda}

\def\intr{\rm int}

\newcommand{\rsquare}

\def\lr{\lambda_{(r)}}

\def\lr1{\hat\lambda_{(r)}}
\def\lm1{\hat\lambda_{(m-r)}}
\def\bee{\begin{equation}}
\def\eee{\end{equation}}

\def\proof{\noindent {\bf Proof}. \ }
\numberwithin{equation}{section}
\title[Exterior Distance Function]{Exterior Distance Function}

\author[R. Polyak]{Roman A. Polyak}

\address[R. Polyak]{ Department of Mathematics\\
The Technion - Israel Institute of Technology\\
 32000 Haifa, Israel}
\email{{\tt rpolyak@techunix.technion.ac.il\and rpolyak@gmu.edu}}



\begin{document}

\begin{abstract}
We introduce and study exterior distance function (EDF) and correspondent exterior point method (EPM) for convex optimization. 

The EDF is a classical Lagrangian for an equivalent problem obtained from  the initial one  by monotone transformation of both the objective function and the constraints. 

The constraints transformation is scaled by a positive scaling parameter. Thus, the EDF  is a particular realization of the Nonlinear Rescaling (NR) principle. 

 Along with the "center", the EDF  has two extra tools: the barrier (scaling) parameter and the vector of Lagrange multipliers.

We show that  EPM  generates  primal - dual sequence, which  converges to the primal - dual solution in value under minimum assumption on the input data. Moreover, the convergence is taking place under any fixed interior point as a "center" and any fixed positive  scaling parameter, just  due to the Lagrange multipliers update. 

If the second order sufficient optimality condition  is satisfied, then the EPM converges with Q-linear rate under any fixed interior point as a  "center" and any fixed, but large enough positive scaling parameter. 
\end{abstract}
\maketitle

\section{Introduction}
The Interior Distance functions (IDFs) were introduced and the Interior Center Methods (ICMs) were  developed  by P.\ Huard in the mid - $60^{th}$  (see \cite{bui_trong66},\cite{huard67},\cite{Huard67}). 

  Later IDFs and correspondent  ICMs were incorporated into SUMT and studied by A.\
Fiacco and G.\ McCormick in \cite{fiacco68}  and other authors (see, for example,  \cite{grossman81},\cite{polyak71} and references therein).

At each step ICM   finds  a central (in a sense)
point of the Relaxation Feasible Set (RFS) and updates the level set 
using the new objective function value. The RFS is the
intersection of the feasible set and the relaxation (level) set of
the objective function at the attained level.  

The ``center" is
sought as a minimizer of the IDF. It is a point in the RFS "most distant" from both the boundary of the objective function level set and the  active constraints.

Interest in  IDFs and correspondent center methods has grown
dramatically after N.\ Karmarkar  published is 1984 his projective
scaling method (see \cite{karmarkar84}).  In fact, his potential function is an IDF and his
method is a Center Method, which generates centers of spheres 
belonging to the interior of the polytope.

Mainly for this reason  the concept of centers became extremely popular in the 80s. Centering and reducing the cost are two basic ideas behind  the Interior Point Methods  (IPMs), which was  the main stream in Modern Optimization for a long time. Centering means to stay away from the boundary of the RFS. An answer to the  basic question: how far from the boundary one should stay in case of LP  was given by Sonnevend in \cite{sonnevend85} (see also \cite{jarre88}) through the definition of  analytic center  of a polytop.  The central path is a curve formed  by  analytic centers. The curve plays an important role in the IPMs (see \cite{gonzaga92}). 

Following the central path J. Renegar in \cite{renegar88} obtained the first  path-following algorithm with $O(\sqrt{n}L)$ number of iterations, versus $O(nL)$ iterations  for the N. Karmarkar's method. 

Soon after C.Gonzaga \cite{gonzaga88} and P. Vaidya \cite{vaidya87} developed algorithms for LP,  based on the centering ideas, with overall complexity $O(n^3L)$ arithmetic operations, which is the best knows result so far.

After Yu. Nesterov and A. Nemirovsky developed their self-concondance theory it becomes evident that path-following methods with polynomial complexity for convex optimization problems is possible  if the RFS can be equipped with self-concordant barrier (see \cite{nesterov94},\cite{nesterov04}). 

If it is not the case, then one can use the classical IDF and correspondent ICM.

The classical  IDF, however, has well known drawbacks: (1) the IDF, its  gradient and Hessian does not exist at the primal solution; (2) the IDF, as well as, the condition number of IDF's Hessian unboundedly grows when the primal approximation approaches the solution. The singularity of the IDF at the solution leads to numerical instability, in particular, in the final phase. 
It means that from some point on, finding an accurate approximation for the  IDF's minimizer is practically an impossible task.

In spite of a long history of IDF and correspondent ICM the fundamental question still is:  \emph{how
 the main idea of center methods: to stay away from
the boundary  consistent with the main purpose of constrained optimization:  
finding a solution on the boundary}.

The issue was partially addressed in \cite{polyak97_2}, where the Modified Interior Distance Functions was introduced and correspondent theory and methods were developed.  The results in \cite{polyak97_2}, however, were obtained  only under the second order sufficient optimality condition.

In this paper we address the issue by introducing the Exterior  Distance Function (EDF) and correspondent Exterior Point Method (EPM). The EDF is a classical Lagrangian for a convex optimization problem  equivalent to the initial one and obtained from the latter by transforming both the objective function and the constraints.

So, EDF is a particular realization of the Nonlinear Rescaling (NR) principle, but the main  EDF results do not follow from NR theory (see \cite{polyak92},\cite{polyak97},\cite{polyak16}).

 We  obtained the basic convergence results under minimum assumptions on the input data.

In contrast to the classical  IDF, the EDF, its gradient and Hessian are  defined on an extended feasible set. It eliminates the singularity of the EDF and its derivatives at the solution.

 The EDF   has two extra tools, which control the  computational process: the positive barrier (scaling) parameter and the vector of Lagrange multipliers.

The EPM  alternates finding the EDF primal minimizer  with Lagrange multipliers update, while both the "center" and the barrier parameter can be fixed or updated from step to step. 

Under a fixed "center"  the EDF resemble the Modified barrier Function (MBF) (see \cite{polyak92}), but  for a problem equivalent to the initial one. Due to the "center" it provides, on the top of the MBF qualities,  an extra one.  By changing the "center" from step to step it is possible  strengthening convergence  results typical for MBF without much extra computational work.

Convergence due to the Lagrange multipliers update allows  keeping  the condition number of the EDF's Hessian stable, which is critical for numerical stability.  This is a fundamental departure from classical IDF theory and methods.

 Under  standard second order sufficient optimality condition  EPM converges with Q-linear rate even  when both  the "center" and the scaling  parameter are fixed, but the parameter  is large enough. Therefore by changing the scaling parameter and/or the "center" from step to step one gets superlinear convergence rate versus sublinear, which is typical for the Classical IDF.

Also for a fixed, but large enough scaling parameter and any fixed interior point as a "center"  the EDF is strongly convex in the neighborhood of  the primal minimizer no matter the objective function and the active constraints  are convex or not.

\section{Problem formulation and basic assumptions}\label{MIDF}

Let  $f:\Rn\to \RR$ be convex and all  $c_i:\Rn\to \RR$, $i=1,...,m$ are concave.

We consider the following convex optimization problem
\begin{equation}\label{Eq1_4.9.16}
f(x^*)=\min\{ f(x)|x\in \Omega\},
\end{equation}
where 
\begin{equation}\label{Eq2_4.9.16}
\Omega=\{x\in\Rn:c_i(x)\geq 0, i=1,...,m\}
\end{equation}
is a feasible set. 

We assume:
\begin{enumerate}
\item[A.] The primal solution set  $X^*=\{x\in \Omega: f(x)=f(x^*)\}$ is not empty and bounded;
\item[B.] Slater condition   
$$\exists x_0\in \Omega: c_i(x_0)>0, \; i=1,...,m$$ holds.
\end{enumerate}
%---------------------------------------------
%\section{Problem reformulation }
Let $y\in\intr\;\Omega$, then the relaxation feasible set (RFS) 
$$\Omega(y)=\{x\in \Omega:f(x)<f(y)\}$$

 is convex  and  bounded for any   given  $y\in \intr\;\Omega$.  It follows from A, convexity $f$, concavity $c_i$, $i=1,...,m$ and Corollary 20 (see  \cite{fiacco68}). 

Also, without losing  generality, we can assume that $ f(x) \geq 0$, because, otherwise, we can replace
$f(x)$ by an equivalent objective function $f(x) := \ln (e^{f(x)} +1)\geq 0$. 

Throughout the  paper  we will use the following well known fact.
\begin{lemma}[Debreu]
Let $A=A^T:\Rn\to\Rn$, $C:\Rn\to\RR^r (n>r)$, $\rank C=r$ and 
$$(Ax,x)\geq \mu(x,x),\;\; \mu>0,\;\;\forall x:Cx=0,$$
then there is $0<\rho<\mu$ and large enough $k_0>0$ that for any $k\geq k_0$ the following inequality
$$((A+kC^TC)x,x)\geq\rho(x,x),\forall x\in\Rn$$
holds true.
\end{lemma}
In the following section we recall some  IDF properties.
%-----------------------------------------------
\section{Classical Interior Distance Function}
Let $y_0 \in$ int $\Omega$ and $\tau_0 = f(y_0)$, then $$ \Omega(\tau_0) =
\{x\in\Omega: f(x) \leq \tau_0\}$$ is the RFS at the level $\tau_0=f(x_0)$.

Let $\tau\in T = \{\tau:
\tau_0 > \tau > \tau^* = f(x^*)$, then Huard's  IDF $H: \Omega (\tau) \times T \to \RR^1$ is
defined by 
\begin{equation}\label{Eq_19.11.16}
H(x, \tau) = - m \ln (\tau-f(x)) - \sum_{i=1}^m \ln c_i(x). \end{equation}

We assume  $\ln t = - \infty$, for $t\leq 0$, then Interior Center Method (ICM) step  finds the ``center" 
\begin{equation}\label{Eq3_19.11.16}
\hat{x} = \hat{x}(\tau) = \argmin\{H(x, \tau)/ x \in \RR^n\}\end{equation}
and replaces $\Omega(\tau)$ by $\Omega(\hat{\tau})$, where 
$\hat{\tau} = f(\hat{x})$. 

From the statement $x \to
\partial \Omega (\tau) \Rightarrow H(x, \tau) \tinf$ follows 
$\hat{x} (\tau) \in$ int $\Omega (\tau) $ for any
$\tau \in T$.  

If the RFS can be equipped with   self - concordant (SC) barrier, that is $H(x,\tau)$ is a SC-function in $x \in   \intr \; \Omega(\tau)$, then following the  central trajectory $\{\hat{x}(\tau), \tau\in T\}$ one  gets an  $\varepsilon$-approximation of   $f(x^*)$ in $O(\sqrt{n}\ln\varepsilon^{-1})$ steps by alternating Newton's step applied for minimization $H(x,\tau)$ with $\tau$ update  (see \cite{nesterov94}, \cite{nesterov04}, \cite{renegar88} and references therein).

If the RFS can not be equipped with  SC barrier, then Classical Huard's IDF and correspondent  ICM  \eqref{Eq3_19.11.16} is used.

  The IDF $F(x,\tau) $, as well as, its gradient and Hessian are singular at $x^*$. 
  
  For any given $\tau\in T$ we have $\lim_{x\to x^{*}}H(x, \tau)=\infty$ and the condition number of the Hessian $\nabla_{xx}H(x,\tau)$ unboundedly grows when $\hat{x}(\tau)\to x^*$, which makes  finding a  reasonable approximation for  primal minimizer, from  some point on, practically  impossible.

Although approximations for the Lagrange multipliers can be found as
a by-product of  ICM, they cannot be effectively used in the
computational process. 

Let us consider the issues briefly. From the
boundedness of  RFS $\Omega (\tau)$ and $\lim_{x\to \partial \Omega(\tau)} H(x,\tau)=\infty$   the unconstrained
minimizer  always exists and $\hat{x} =\hat{x}(\tau)\in$ int $\Omega (\tau)$. So we have 
\begin{equation}\nabla_x H(\hat{x},
\tau) = \frac{m}{\tau - f(\hat{x})} \ \nabla f(\hat{x}) - \sum_{i=1}^m \frac{\nabla
c_i(\hat{x})}{c_i(\hat{x})} = 0 \label{MIDF_eq2.1}
 \end{equation} or
 \begin{equation} \nabla f(\hat{x}) -
\sum_{i=1}^m \frac{\tau - f(\hat{x})}{m \, c_i(\hat{x})} \ \nabla c_i(\hat{x}) = 0
\,. \label{MIDF_eq2.2} 
\end{equation}

Let   \begin{equation} \hat{\l}=(\hat{\l}_i =
\hat{\l}_i (\tau) = (\tau- f(\hat{x})) (m \, c_i (\hat{x}))^{-1},
\quad i=1,\ldots, m )\label{MIDF_eq2.123} \end{equation} be the vector of  Lagrange multipliers. 

Vector $\hat{\l}$ is positive 
because $f(\hat{x}) < \tau$ and all $c_i
(\hat{x}) > 0$. %Therefore, the vector $\hat{\l} \equiv \hat{\l} (\tau) = (\hat{\l}_i
%(\tau)$, $i=1,\ldots, q)$ is positive.
The systems \eqref{MIDF_eq2.1} and (\ref{MIDF_eq2.2}) can be rewritten as follows:

\begin{equation}
\nabla_x H(x,\tau)= \frac{m}{\tau-f(\hat{x})}\nabla_x L(\hat{x}, \hat{\l}) = \frac{m}{\tau-f(\hat{x})}(\nabla f(\hat{x})-\nabla c(\hat{x})^T\hat{\l})=0 , \end{equation}
 where  $\nabla c(x) = J(c(x))$ is the $m\times n$ Jacobian of $c(x)=(c_1(x),...,c_m(x))^T$. 
 
 From (\ref{MIDF_eq2.123}) we have
\begin{equation}
\label{MIDF_eq5}
\hat{\l}_i c_i (\hat{x}) = (\tau - f(\hat{x})) m^{-1}, \ i=1,\ldots, m \, . \end{equation}
Summing up \eqref{MIDF_eq5}, we obtain $$ \hat{\l}^T c(\hat{x}) =
\sum_{i=1}^m \hat{\l}_i c_i(\hat{x}) = \tau - f(\hat{x}) \, . $$ 

From \eqref{MIDF_eq2.1} and  $\tau> f(\hat{x})$ follows  $\tau - f(\hat{x}) \tzer,$ when $ \tau\to\tau^{*}$, because $f$ is bounded from below.

Vector $\hat{x}\in\mbox{ int } \Omega$ is primal feasible,   vector  $\hat{\l}\in R^m_{++},$ is dual feasible and from \eqref{MIDF_eq5} follows asymptotic complementarity condition $$\lim_{\tau\to\tau^*}\hat{\l}_i(\tau)c_i(\hat{x}(\tau))\to 0,\;i=1,...,m.$$

 To simplify considerations  we assume at this point that the second order sufficient  optimality  condition for the problem (\ref{Eq1_4.9.16}) is satisfied. Then   the primal-dual solution $(x^*, \l^*)$ is unique.
Therefore,
$$
 \quad \lim_{\tau \to \tau^*} \hat{x}
(\tau) = x^*, \quad  \lim_{\tau \to \tau^*} \hat{\l} (\tau) = \l^* .$$
Let $ I^*=\{i:c_i(x^*)=0\}=\{1,...,r\}$ be the active constraints set.

For  Hessian $\nabla^2_{xx} H(x,\tau)$ at $x=\hat{x}$ we obtain 
 \begin{eqnarray*} \nabla_{xx}^2 H(x,
\tau)_{/x=\hat{x}} & = & m(\tau-f(\hat{x}))^{-1}\left[(\tau-f(\hat{x}))^{-1} \nabla f (\hat{x})
\nabla f^T(\hat{x}) + \nabla^2 f(\hat{x}) \right . \\
&& \left . - \sum_{i=1}^m \frac{(\tau-f(\hat{x}))}{m} \frac{\nabla^2
c_i(\hat{x})}{c_i(\hat{x})} + \sum_{i=1}^m \frac{(\tau-f(\hat{x}))}{mc_i^2(\hat{x})}
\nabla c_i (\hat{x}) \nabla c_i (\hat{x})^T \right] \\
& = & m(\tau-f(\hat{x}))^{-1} \lsb \nabla_{xx}^2 L(\hat{x}, \hat{\l}) + \nabla c
(\hat{x})^T C^{-1} (\hat{x}) \hat{\Lambda} (\tau) \nabla c(\hat{x}) \right . \\
&& + \left . \ (\tau-f(\hat{x}))^{-1} \nabla f(\hat{x}) \nabla
f(\hat{x})^T \right],
\end{eqnarray*}
where  $C(x) = [\diag \, c_i(x)]_{i=1}^m$ and $ \Lambda (\tau) =
[diag \, \l_i (\tau) ]_{i=1}^m $ are diagonal matrices and $L(x,\l)=f(x)-\sum^m_{i=1}\l_ic_i(x)$ is the  Lagrangian for  problem (\ref{Eq1_4.9.16}).

In view of $\hat{x} = \hat{x} (\tau) \to x^*$ and $\hat{\l} = \hat{\l}
(\tau) \to \l^*$   for $\tau$ close to $\tau^*$ we have
\begin{eqnarray*}
\nabla_{xx}^2 H(\hat{x}, \tau) & \approx & m(\tau - f(x^*))^{-1} \lsb \nabla_{xx}^2
L(x^*,
\l^*) + \nabla c(x^*)^T \Lambda^*  C^{-1} (\hat{x}) \nabla c(x^*) \right . \\
&& + \ \left . (\tau-f(x^*))^{-1} \nabla f (x^*) \nabla f(x^*)^T \rsb .
\end{eqnarray*}

From the K-K-T condition  
$$
\nabla f(x^*) = \sum_{i=1}^r \l_i^* \nabla c_i(x^*), $$ follows 
$$  \forall \, u: \nabla c_{(r)} (x^*) \, u=0 \Rightarrow
(\nabla f(x^*),u) = 0 \,  , $$
where $ \nabla  c_{(r)}(x^*)=J(c_{(r)}(x^*))$ is $r\times n$ Jacobian of the vector-function $c_{(r)}(x)=(c_1(x),...,c_r(x))^T$, which corresponds to the active constraints. 

 Hence, for
$\forall \, u: \nabla c_{(r)} (x^*)u = 0$, we obtain
$$ \lb \nabla_{xx}^2  H(\hat{x}, \tau)u,u \rb \approx m(\tau - f(x^*))^{-1} \, \times
$$
\begin{equation}
 \lb \lb \nabla_{xx}^2 L(x^*, \l^*) + \nabla c_{(r)}^T (x^*) \Lambda_{(r)}^*
C_{(r)}^{-1} (\hat{x}) \nabla c_{(r)} (x^*) \rb u,u \rb, \label{MIDF_eq2.123a}
\end{equation}
where
$$
\Lambda_{(r)}^* = \diag \, (\l_i^*)_{i=1}^r, \ C_{(r)} (x) = \diag
\, (c_i(\hat{x}))_{i=1}^r \, . $$ From the second order sufficient  condition follows $\l^{*}_i>0$, $i=1,...,r$, also from $\hat{x} = \hat{x} (\tau)
\to x^*$ follows  $c_i (\hat{x}) \tzer$, \ $i=1,\ldots, r$.  

Therefore 
\begin{equation}
\lim_{\tau \to f(x^*)} M_i (\tau) = \lim_{\tau \to f(x^*)}
\hat{\l}_i(\tau) c_i^{-1} (\hat{x}(\tau)) = \infty \, , \quad i=1, \ldots, r \,
. \label{MIDF_eq2.124} 
\end{equation}
For  $\tau_0>\tau>\tau^*=f(x^*)$ close to $\tau^*$  from
 Debreu's lemma  with $A=\nabla^2_{xx}L(x^{*},\l^*)$ and $C=\Lambda_{(r)}^{*\frac{1}{2}}(C_{(r)}(\hat{x}))^{-\frac{1}{2}}\nabla c_{(r)}(x^*)$  follows existence  of  $\rho
> 0$, such that  
$$
\mu(\tau)=\mbox{mineigenval } \nabla_{xx}^2 H(\hat{x}, \tau) \to \rho
$$ 
when $\tau\to\tau^*$.
 
On the other hand, from \eqref{MIDF_eq2.124} follows 
$$
M(\tau)=\mbox{maxeigenval } \nabla_{xx}^2 H(\hat{x}, \tau) \to\infty, $$ 
 when $\tau\to\tau^*$.

 %and the condition number
Therefore 
$$\mbox{cond}\nabla^2_{xx} H(\hat{x},\tau)=M(\tau)\mu^{-1}(\tau)\rightarrow \infty$$
 when $\tau \to \tau^*$.

The ill-conditioning of the Hessian $\nabla_{xx}^2 H(\hat{x}, \tau)$ is much more
critical in nonlinear optimization than in LP.  In case of LP, the term
$\nabla_{xx}^2 L(x, \l)$ in the expression of the Hessian  $\nabla_{xx}^2 H(\hat{x},
\tau)$ disappears and by rescaling the input data properly, one can, to some extent,
eliminate the ill-conditioning effect.

In nonlinear optimization, the situation is completely different and the
ill-conditioning is an important issue, in particular, when  solution with high accuracy is required.

In the  following section we introduce and study  the EDF, which eliminates the basic drawbacks of the Classical  IDF.

%-----------------------------------------------------------
\section{Exterior Distance Function}
  For a given $y\in\intr\;\Omega$ let us consider the following problem 
  \begin{equation}\label{Eq3_4.9.16}
  F(x^*,y)=\min\{F(x,y)|c_i(x)\geq 0, i=1,...,m\}, 
  \end{equation}
 where
 $$ F(x,y)=-\ln\Delta(x,y)=-\ln(f(y)-f(x)).$$
 
 For any   $y\in\intr\;\Omega$ the function $F$ is convex  and monotone decreasing together with $f$ for  $x\in\Omega(y)=\{x\in\Omega: f(x)\leq f(y)\}$, therefore the solution $x^*\in\Omega(y)$ of the problem (\ref{Eq3_4.9.16}) belongs to $X^*$ and vice versa any $x^*\in X ^*$ solves (\ref{Eq3_4.9.16}), that is problems \eqref{Eq1_4.9.16} and    (\ref{Eq3_4.9.16}) are equivalent. 
 
 In what is following  we consider  the problem (\ref{Eq3_4.9.16}) instead  (\ref{Eq1_4.9.16}) .
 
 The correspondent to  (\ref{Eq3_4.9.16}) Lagrangian $L_y:\Rn\times \Rm_+\to\RR$  is given by 
 \begin{equation}\label{Eq4_4.9.16}
 L_y(x,\l)=F(x,y)-\sum^m_{i=1}\l_ic_i(x).
 \end{equation}
 
 The correspondent to (\ref{Eq3_4.9.16}) dual function $d_y:\RR^m_+\to \RR$ is  $$d_y(\l)=\inf_{x\in\RR} L_y(x,\l)$$  and 
 \begin{equation}\label{Eq5_4.9.16}
 d_y(\l^*)=\max\{d_y(\l)|\l\in\Rm_+\}
 \end{equation}
is the dual to \eqref{Eq3_4.9.16} problem.
 
 Let $\psi(t)=\ln(t+1)$ and $k>0$, then the original set of constraints $$c_i(x)\geq 0,\quad i=1,...,m$$ is equivalent to the following set 
 \begin{equation}\label{Eq6_4.9.16}
 k^{-1}\psi(kc_i(x))=k^{-1}\ln(k c_i(x)+1)\geq0,\; i=1,...,m,
 \end{equation}
 therefore for any given $k>0$ we have
 $$\Omega=\{x\in\Rn:k^{-1}\psi(k c_i(x))=k^{-1}\ln(kc_i(x)+1)\geq 0, \; i=1,...,m\},$$
 
 and for any given $y\in \intr\;\Omega$ and $k>0$  the  problem 
 \begin{equation}\label{Eq7_4.9.16}
 F(x^*,y)=\min\{F(x,y)|k^{-1}\ln(kc_i(x)+1)\geq0,\; i=1,..,m\}
 \end{equation}
 is equivalent to  (\ref{Eq3_4.9.16}).
 
 Let us fix $y\in \intr\;\Omega$, then the following  extension
 \begin{equation}\label{Eq8_4.9.16}
 \Omega_{-k^{-1}}(y)=\{x\in\Rn:c_i(x)\geq-k^{-1},i=1,...,m,\; f(y)>f(x)\}
 \end{equation}
 of $\Omega(y)$ is convex and bounded due to convexity $f$, concavity $c_i$, $i=1,...,m$ boundedness $\Omega(y)$ and Corollary 20  (\cite{fiacco68})  and so is  the following contruction of $\Omega(y)$
 \begin{equation}\label{Eq9_4.9.16}
 \Omega_{\gamma}(y)=\{x\in\RR:c_i(x)\geq \gamma,\;i=1,..,m, \;\;f(y)>f(x)\}.
 \end{equation}
  The set $\Omega_{\gamma}(y)$ for small $\gamma>0$ is not empty due to the Slater condition.
 
 Let us fix $ y\in \intr\;\Omega$  and $k>0$,  then Lagrangian $\mathcal{L}_y:\Rn\times\Rm_+\times\RR_{++}\to\RR$  
 \begin{equation}\label{Eq10_4.9.16}
 \mathcal{L}_y(x,\l,k)=F(x,y)-k^{-1}\sum^m_{i=1}\l_i\ln(kc_i(x)+1),
 \end{equation}
for   problem \eqref{Eq7_4.9.16} we call  the exterior distance function (EDF).
 
Thus, EDF is a particular realization of the  NR principle (see \cite{polyak92}, \cite{polyak97}, \cite{polyak16} and references therein).
 
 Let us consider the second order sufficient optimality conditions for  problem (\ref{Eq3_4.9.16}).
 
  There exists  $\mu>0$, such that 
 \begin{equation}\label{Eq11_4.9.16}
 (\nabla^2_{xx}L_y(x^*,\l^*)u,u)\geq\mu(u,u),\forall u:\nabla c_{(r)}(x^*)u=0
 \end{equation}
 and 
 \begin{equation}\label{Eq12_4.9.16}
 \rank \nabla c_{(r)}(x^*)=r.
 \end{equation}

 We conclude the section by pointing out some EDF properties at the KKT's point $(x^*,\l^*)$.
 
 First of all, $\mathcal{L}$ is convex in $x\in \Omega_{-k^{-1}}(y)$ for any given  $\intr\;\Omega$, $k>0$ and $\l\in \Rm_+$. 
 \begin{proposition}
 For a given $y\in\intr\;\Omega$ and $k>0$  and any  KKT's point $(x^*,\l^*)$ we have:
 \begin{enumerate}
 \item[$1^0$]$\mathcal{L}_y(x^*,\l^*,k)=F(x^*,y)=-\ln(f(y)-f(x^*))$
 
 or
 $$f(x^*)=f(y)-e^{-F(x^*,y)};$$
 \item[$2^0$] $\nabla_x\mathcal{L}_y(x^*,\l^*,k)=\Delta^{-1}(x^*,y)\nabla f(x^*)-\sum^m_{i=1}(kc_i(x^*)+1)^{-1}\l^*_i\nabla c_i(x^*)=$
 $$\Delta^{-1}(x^*,y)\nabla f(x^*)-\sum^m_{i=1}\l^*_i\nabla c_i(x^*)=\nabla_xL_y(x^*,\l^*)=0;$$
 \item[$3^0$] $\nabla^2_{xx}\mathcal{L}_y(x^*,\l^*,k)=\Delta^{-2}(x^*,y)\nabla f(x^*)\nabla f^T(x^*)+\Delta^{-1}(x^*,y)\nabla^2 f(x^*)$
 $$-\sum^m_{i=1}\l^*_i\nabla^2c_i(x^*)+k\nabla c(x^*)^T\Lambda^*\nabla c(x^*)=\nabla^2_{xx}L_y(x^*,\l^*)+k\nabla c(x^*)^T\Lambda^*\nabla c(x^*)=$$
 $$\nabla^2_{xx}L_y(x^*,\l^*)+k\nabla c_{(r)}(x^*)^T\Lambda^*_{(r)}\nabla c_{(r)}(x^*),$$
 \end{enumerate}
  where  $\Lambda^*_{(r)}=\diag (\l^*_i)^r_{i=1}$, $\l^*_i=0$, $i=r+1,...,m$.
  \end{proposition}
  Properties  $1^0-3^0$ follow from the definition of EDF (\ref{Eq10_4.9.16}) and complementarity condition 
   \begin{equation}\label{Eq13_4.9.16}
   \l^*_ic_i(x^*)=0, i=1,...,m.
   \end{equation}
  The fundamental difference between EDF  \eqref{Eq10_4.9.16} and  the   Huard's IDF (\ref{Eq_19.11.16}) follows from  $1^0-3^0$.

 First, the $\mathcal{L}_y(x,\l,k)$ is defined at  the solution together with its gradient  and Hessian.
 
 Second, from $2^0$ follows that for any given $y \in\intr\;\Omega$ and $k>0$ the optimal solution of (\ref{Eq3_4.9.16}) can be found by solving one smooth unconstrained optimization problem
 \begin{equation}\label{Eq14_4.9.16}
  \min_{x\in\Rn} \mathcal{L}_y(x,\l^*,k)=\mathcal{L}_y(x^*,\l^*,k).
 \end{equation}
 It means that $\mathcal{L}_y(x,\l^*,k)$ is an exact smooth approximation for the following  non-smooth problem 
 \begin{equation}\label{Eq15_4.9.16}
  \min_{x\in\Rn}\max\{F(x,y)-F(x^*,y),-c_i(x),\;i=1,...,m\},
 \end{equation}
 which is for any given $y\in \intr \;\Omega$ is equivalent to (\ref{Eq3_4.9.16}).
 
 Third, from $3^0$  for any $u\in\Rn$ follows 
 $$(\nabla^2_{xx}\mathcal{L}_y(x^*,\l^*,k)u,u)=((\nabla^2_{xx}L_y(x^*,\l^*)+k\nabla c_{(r)}(x^*)^T\Lambda^*_{(r)}\nabla c_{(r)}(x^*))u,u).$$
 \begin{proposition}\label{P2_4.11.16}
 Under the second order sufficient optimality condition (\ref{Eq11_4.9.16})-(\ref{Eq12_4.9.16}),  for any given $y \in\intr\;\Omega$, $k_0>0$ large enough and any $k\geq k_0$ there exists $0<\rho <\mu$ such that 
  $$(\nabla^2_{xx}\mathcal{L}_y(x^*,\l^*,k)u,u)\geq \rho(u,u),\; \forall u\in \Rn.$$
 \end{proposition}
 Proposition \ref{P2_4.11.16} follows from the  second order sufficient optimality condition (\ref{Eq11_4.9.16})-(\ref{Eq12_4.9.16}) and  Debreu's Lemma with $$A=\nabla^2_{xx}L_y(x^*,\l^*) \mbox{ , }  C=\Lambda^{*\frac{1}{2}}_{(r)}\nabla c_{(r)}(x^*).$$
 
 In other words, for any fixed $y\in\intr\;\Omega$ as a "center" and any $k\geq k_0$ the EDF $\mathcal{L}_y(x,\l^*,k)$ is strongly convex in the neighborhood of $x^*$ no matter if $f$ and $-c_i$, $i=1,...,m$ are convex or not.
 
 The EDF is related to the Classical Huard's interior distance function $H(x,\tau)$ as MBF (see \cite{polyak92}) is to the classical R. Frisch's  barrier function $F(x,k)=f(x)-k^{-1}\sum^m_{i=1}\ln c_i(x)$ (see \cite{frisch55}). 
 
It is worth mentioning that relatively to the   MBF features the EDF has one extra tool-the "center", which we will use later to improve convergence properties.
 
 The EDF properties  lead to a new multipliers method, which converges under any fixed $y \in \intr\;\Omega$ as a "center" and any fixed scaling parameter $k>0$, just due to the Lagrange multipliers update. This is a fundamental departure from the Classical IDF theory (see \cite{bui_trong66},\cite{fiacco68},\cite{grossman81},\cite{huard67},\cite{Huard67},\cite{polyak71}).
 
 %-------------------------------------------------
 \section{Exterior Point Method}
 The EPM  at each  step finds the primal minimizer of $\mathcal{L}_y$ following by Lagrange multipliers updates.
 
We start with   $y\in \intr\;\Omega$ as a given fixed "center", fixed scaling  parameter $k>0$ and  initial Lagrange multipliers  vector $\l_0=e=(1,...,1)^T\in\Rm_{++}$. 

Let  the  primal-dual approximation $(x_s,\l_s)$ has been found already.

 The  approximation $(x_{s+1},\l_{s+1})$ we find  by the following operations 
 
$$ x_{s+1}:\nabla_x\mathcal{L}_y(x_{s+1},\l_s,k)$$
 
 \begin{equation}\label{Eq16_4.9.16}=\Delta^{-1}(x_{s+1},y)\nabla f(x_{s+1})-\sum^m_{i=1}\l_{i,s}\psi^{'}(kc_i(x_{s+1}))\nabla c_i(x_{s+1})=0
 \end{equation}
 \begin{equation}\label{Eq17_4.9.16}
 \l_{s+1}:\l_{i,s+1}=\l_{i,s}\psi^{'}(kc_i(x_{s+1}))=\l_{i,s}(kc_i(x_{s+1})+1)^{-1},\; i=1,...,m.
 \end{equation}
 The key ingredient  of the EPM (\ref{Eq16_4.9.16})-(\ref{Eq17_4.9.16})  convergence analysis is its  equivalence to the proximal point method with $\varphi$-divergence distance function for the dual problem
 \begin{theorem}\label{T1_5_9.16}
 If condition A and B hold, $f$, $c_i\in C^1$, $i=1,...,m$, $f$ is convex and all $c_i$, $i=1,...,m$ are concave,  then EPM (\ref{Eq16_4.9.16})-(\ref{Eq17_4.9.16}) is:
 \begin{enumerate}
 \item[1)] well defined;
 \item[2)] equivalent to the following proximal point method
 \begin{equation}\label{Eq18_4.9.16}
 d_y(\hat{\l})-k^{-1}D(\hat{\l},\l)=\max\{d_y(u)-k^{-1} D(u,\l)|u\in \Rm_+\},
 \end{equation}
 where 
 $$D(u,\l)=\sum^m_{i=1}\l_i\varphi(u_i/\l_i)$$
 is $\varphi$-divergence distance function based on the kernel $\varphi=-\psi^*$, where $\psi^*$ is Legendre transform of $\psi$
 \end{enumerate}
 \end{theorem}
\proof 
\begin{enumerate}
\item[1)] Due to convexity $f$, concavity $c_i$, $i=1,...,m$, Slater condition,  boundedness $\Omega(y)$ and properties of log-barrier function the recession cone of $\Omega(y)$ is empty, that is we have 
$$\lim_{t\to\infty}\mathcal{L}_y(x+td,\l,k)=\infty$$
 for any  $ d\neq 0 $ from $\Rn$, $y\in\intr\;\Omega$, $k>0$ and $\l\in\Rm_{++}$.
 
Hence, there exists $x_{s+1}\in\Rn:$
 $$\mathcal{L}_y(x_{s+1},\l_s,k)=\min\{\mathcal{L}_y(x,\l_s,k)|x\in\Rn\},$$
 thus (\ref{Eq16_4.9.16}) holds.
 
 From $\ln t=-\infty$ for $t\leq 0$  and (\ref{Eq16_4.9.16}) follows  $kc_i(x_{s+1})+1>0$, $i=1,...,m$, therefore from (\ref{Eq17_4.9.16})  we have
 $$\l_s\in\Rm_{++}\Rightarrow \l_{s+1}\in\Rm_{++}.$$
 Hence, method (\ref{Eq16_4.9.16})-(\ref{Eq17_4.9.16}) is well defined.
 
 \item[2)] From (\ref{Eq16_4.9.16}) and (\ref{Eq17_4.9.16}) follows
 $$\nabla_x\mathcal{L}_y(x_{s+1},\l_s,k)=\Delta^{-1}(x_{s+1},y)\nabla f(x_{s+1})-\sum^m_{i=1}\l_{i,s+1}\nabla c_i(x_{s+1})=$$
 $$\nabla_x L_y(x_{s+1},\l_{s+1})=0.$$
 Therefore 
 $$\min_{x\in\Rn}L_y(x,\l_{s+1})=L_y(x_{s+1},\l_{s+1})=d_y(\l_{s+1}).$$
 The subdifferential $\partial d_y(\l_{s+1})$ contains $-c(x_{s+1})$, that is 
 \begin{equation}\label{Eq19_4.9.16}
 0\in c(x_{s+1})+\partial d_y(\l_{s+1}).
 \end{equation}
 From (\ref{Eq17_4.9.16}) we have
 $$\psi^{'}(k c_i(x_{s+1}))=\l_{i,s+1}/\l_{i,s},\;i=1,...,m.$$
 Also $\psi^{''}(kc_i(x_{s+1}))\neq 0$, therefore the inverse function $\psi^{'-1}$ exists and 
 \begin{equation}\label{Eq20_4.9.16}
 c_i(x_{s+1})=k^{-1}\psi^{'-1}(\l_{i,s+1}/\l_{i,s}).
 \end{equation}
 From (\ref{Eq20_4.9.16})  and Legendre identity $\psi^{'-1}\equiv \psi ^{*'}$ follows 
 \begin{equation}\label{Eq21_4.9.16}
 c_i(x_{s+1})=k^{-1}\psi^{*'}(\l_{i,s+1}/\l_{i,s}),\; i=1,...,m.
 \end{equation}
 From (\ref{Eq19_4.9.16}) and (\ref{Eq21_4.9.16}) we obtain 
 \begin{equation}\label{Eq22_4.9.16}
 0\in \partial d_y(\l_{s+1})+k^{-1}\sum^m_{i=1}\psi^{*'}(\l_{i,s+1}/\l_{i,s})e_i,
 \end{equation}
 where $e_i=(0,...,1,...,0)^T\in\Rm_+$.
 \end{enumerate}
 Let $\varphi=-\psi^{*}$, then (\ref{Eq22_4.9.16}) we can rewrite as follows
 $$0\in \partial d_y (\l_{s+1})-k^{-1}\sum^m_{i=1}\varphi^{'}(\l_{i,s+1}/\l_{i,s})e_i,$$
 which is the optimality condition for $\l_{s+1}$ to be the solution in (\ref{Eq18_4.9.16}) with $\l=\l_s$, that is
 $$d_y(\l_{s+1})-k^{-1}\sum^m_{i=1}\l_{i,s}\varphi(\l_{i,s+1}/\l_{i,s})=$$
 \begin{equation}\label{Eq23_4.9.16}
 \max\{d_y(u)-k^{-1}\sum^m_{i=1}\l_{i,s}\varphi(u_i/\l_{i,s})|u\in\Rm_{++}\}
 \end{equation}
 $$=\max\{d_y(u)-k^{-1}D(u,\l_s)|u\in\Rm_{++}\},$$
 where $\psi^{*}(s)=\ln s-s+1$. 
 
 Therefore 
for the  kernel  $\varphi(s)=-\psi^*(s)$ of the  $\varphi$- divergence distance function
$$D(u,\l)=\sum^m_{i=1}\l_i\varphi(u_i/\l_i)=\sum^m_{i=1}[-\l_i\ln u_i/\l_i+u_i-\l_i]$$

 we have:

\begin{equation}
\begin{split}
& a) \;\;\;\varphi(s)=-\ln s+s-1\geq 0,\;\;\;\forall s>0\\
& b)  \;\;\;\; \min_{s>0}\varphi(s)=\varphi(1)=\varphi^{'}(1)=0.
\end{split}
\end{equation}
 In fact, $D(u,\l)$  is the   Kullback-Leibler $\varphi$- divergence distance function (see, for example,  \cite{polyak16}). The proof of Theorem \ref{T1_5_9.16} is completed.
 
 Let $X$ and $Y$ be two bounded and closed sets in $\Rn$ and $d(x,y)=\|x-y\|$ is Euclidean distance between $x\in X$ and $y\in Y$. Then the Hausdorff distance between $X$ and $Y$ is defined as follows
 $$d_H(X,Y):=\max\{\max_{x\in X}\min_{y\in Y} d(x,y),\max_{y\in Y}\min_{x\in X}d(x,y)\}=$$
 $$\max\{\max_{x\in X} d(x,Y),\max_{y\in Y} d(X,y)\}.$$
  For compact sets  $X$ and $Y$  we have 
  \begin{equation}\label{Eq25_4.9.16}
  d_H(X,Y)=0\Leftrightarrow X=Y.
  \end{equation}
 Let  $Q\subset \Rm_+$ be a compact set, $\hat{Q}=\Rm_+\setminus Q$, $S(u,\varepsilon)=\{v\in \Rm_+:\|u-v\|\leq \varepsilon\}$ and 
 $$\partial Q=\{u\in Q|\exists v\in Q: v\in S(u,\varepsilon),\; \exists \hat{v}\in \hat{Q}:\hat{v}\in  S(u, \varepsilon),\; \forall \varepsilon >0\}$$
 be the boundary of $Q$.
 
 For convex and compact sets $A\subset B \subset C$ the inequality
 \begin{equation}\label{Eq26_4.9.16}
 d_H(A,\partial B)< d_H(A, \partial C)
 \end{equation}
 follows from the definition of Hausdorff distance.
 
 For the dual sequence $\{\l_s\}^{\infty}_{s=0}$ we consider the dual level sets $\Lambda_s=\{\l\in \Rm_+:d(\l)\geq d(\l_s)\},$ which are convex, due to concavity $d$, and bounded, due to the boundedness $L^*$,which is,  in turn, a consequence of Slater's condition. Let $\partial \Lambda_s=\{\l\in\Lambda_s:d(\l)=d(\l_s)\}$ be the boundary of $\Lambda _s$.
 %-------------------------------------------------------------------------
 \section{Convergence of the EPM}
 The following Theorem establishes  convergence of the EPM under minimum assumptions on the input data, just due to the  Lagrange multipliers update.
 \begin{theorem}\label{T_29.11.6}
 Under assumptions of Theorem \ref{T1_5_9.16} for any  fixed  $y\in \intr\;\Omega$, as a  "center", any  scaling parameter $k>0$ and any $\l_0\in \Rm_{++}$ the EPM  \eqref{Eq16_4.9.16}-\eqref{Eq17_4.9.16} generates  primal-dual sequence $\{x_s,\l_s\}^{\infty}_{s+0}$ that:
 \begin{enumerate}
 \item[1)] $d_y(\l_{s+1})>d_y(\l_s),\; s\geq 0$
 \item[2)] $\lim_{s\to\infty}d_y(\l_s)=d_y(\l^*),\; \lim_{s\to\infty} F(y,x_s)=F(y,x^*)$
 \item[3)] $\lim_{s\to\infty} d_H( \partial \Lambda_s,L^*)=0$
 \item[4)] there exists a subsequence $\{s_l\}^{\infty}_{l=1}$ such that  for $\bar{x}_l=\sum^{s_{l+1}}_{s=s_l}(s_{l+1}-s_l)^{-1}x_s$ we have $\lim_{l\to \infty}\bar{x}_l=\bar{x}\in X^*$, that is the primal sequence converges  to the primal solution in the ergodic sence.
 \end{enumerate}
 \end{theorem}
 \proof
  \begin{enumerate}
 \item[1)] 
 From $\varphi(1)=0$ and (\ref{Eq23_4.9.16}) with $u=\l_s$ follows 
 \begin{equation}\label{Eq27_4.9.16}
 d_y(\l_{s+1})\geq d_y(\l_s)+k^{-1}\sum^m_{i=1}\l_{i,s}\varphi(\l_{i,s+1}/\l_{i,s})
 \end{equation}
 $$=d_y(\l_s)+k^{-1}D(\l_{s+1},\l_s).$$
 From $\varphi(t)\geq 0,$ $\forall t>0$, (\ref{Eq27_4.9.16}) and $\l_s\in \Rm_{++}$ follows
 \begin{equation}\label{Eq28_4.9.16}
 d_y(\l_{s+1}) \geq d_y(\l_s).
 \end{equation}
Moreover, $d_y(\l_{s+1})>d_y(\l_s)$  unless $\varphi(\l_{i,s+1}/\l_{i,s})=0$ for all $i=1,...,m$, which leads to $\l_{s+1}=\l_s=\l^*$.
 
 \item[2)]    The monotone increasing sequence
 $\{d_y(\l_s)\}^{\infty}_{s=0}$ is bounded from above by the optimal value of the primal objective function $F(y,x^*)=-\ln (f(y)-f(x^*)),$ 
 therefore there exists $\lim_{s\to \infty} d_y(\l_s)=\bar{d}\leq F(y,x^*).$

 Our next step is to show $\bar{d}=F(y,x^*)$. 
 
 From $-c(x_{s+1})\in \partial d_y(\l_{s+1})$ and concavity of the dual function $d_y$  follows 
 $$d_y(\l)-d_y(\l_{s+1})\leq (-c(x_{s+1}), \l-\l_{s+1}), \; \forall  \l \in \Rm_{++}.$$
 For $\l=\l_s$ we obtain 
 \begin{equation}\label{Eq29_4.9.16}
 d_y(\l_{s+1})-d_y(\l_s)\geq(c(x_{s+1}),\l_s-\l_{s+1}).
 \end{equation}
 From the update formula (\ref{Eq17_4.9.16}) we have 
 \begin{equation}\label{Eq30_4.9.16}
 \l_{i,s}-\l_{i,s+1}=kc_i(x_{s+1})\l_{i,s+1},\; i=1,...,m.
 \end{equation}
 Therefore from (\ref{Eq29_4.9.16}) and (\ref{Eq30_4.9.16}) follows 
 \begin{equation}\label{Eq31_4.9.16}
 d_y(\l_{s+1})-d_y(\l_s)\geq k\sum^m_{i=1}c_i^2(x_{s+1})\l_{i,s+1}.
 \end{equation}
 From boundedness of $L^*$ and concavity $d_y$ follows boundedness of the initial dual level set 
 $$\Lambda_0=\{\l\in\Rm_+:d_y(\l)\geq d_y(\l_0)\}.$$
 From the dual monotonicity (\ref{Eq28_4.9.16}) and boundedness $\Lambda_0$ follows boundedness of the dual sequence $\{\l_s\}^{\infty}_{s=0}\subset \Lambda_0$.
 
 Therefore there exists $L=\max_{i,s} \l_{i,s}$ and from (\ref{Eq31_4.9.16}) follows 
 \begin{equation}\label{Eq32_4.9.16}
 d_y(\l_{s+1})-d_y(\l_s)\geq kL^{-1}(c(x_{s+1}),\l_{s+1})^2.
 \end{equation}
 Summing up (\ref{Eq32_4.9.16}) from $s=1$ to $s=N$ we obtain 
 $$d_y(\l^*)-d_y(\l_0)\geq d_y(\l_{N+1})-d_y(\l_0)\geq kL^{-1}\sum^N_{s+1}(\l_s,c(x_s))^2.$$
 It leads to the asymptotic complementarity condition 
 \begin{equation}\label{Eq33_4.9.16}
 \lim_{s\to\infty}(\l_s,c(x_s))=0.
 \end{equation}
 Summing up (\ref{Eq27_4.9.16})  from $s=0$ to $s=N$ we obtain
 $$d_y(\l^*)-d_y(\l_0)\geq d_y(\l_N)-d_y(\l_0)\geq k^{-1}\sum^N_{s+1}D(\l_{s+1},\l_s),$$
 therefore  $\lim_{s\to\infty}D(\l_{s+1},\l_s)=0.$
 
 The diminishing divergence between two sequential Lagrange multipliers vectors leads us to believe that under any given $y\in\intr\;\Omega$ as a "center" and any given scaling parameter $k>0$ the map
 $$\l\to\hat{\l}(\l,k)=(kC(\hat{x})+I^m)^{-1}\l,$$ has a fixed point $\l^*$, where  $C(\hat{x})=\diag (c_i(\hat{x}))^m_{i=1}$, $I^m$- identical matrix in $\Rm$ and $$\hat{x}(\l,k)\equiv\hat{x}:\nabla_x\mathcal{L}_y(\hat{x},\l,k)=0.$$
 
 First, let us show that for any $\l^*\in\Lambda^*$ the sequence $\{D(\l_s,\l^*)\}^{\infty}_{s=0}$  is monotone decreasing.
 
 We assume $x\ln x=0$ for $x=0$, then 
 $$D(\l_s,\l^*)-D(\l_{s+1},\l^*)=\sum^m_{i=1}(\l^*_i\ln\l_{i,s+1}/\l_{i,s}+\l_{i,s}-\l_{i,s+1}).$$
 Using the update formula (\ref{Eq17_4.9.16})  we obtain 
 \begin{equation}\label{Eq34_4.9.16}
D(\l_s,\l^*)-D(\l_{s+1},\l^*)=
 \end{equation}
 $$\sum^m_{i=1}\l^*_i\ln(kc_i(x_{s+1})+1)^{-1}+k\sum^m_{i=1}\l_{i,s+1}c_i(x_{s+1}).$$
 From $-\ln(1+t)\geq -t,\; \forall t>-1$ and \eqref{Eq34_4.9.16} follows
 \begin{equation}\label{Eq35_4.9.16}
 D(\l_s,\l^*)-D(\l_{s+1},\l^*)\geq 
 \end{equation}
 $$k\sum^m_{i=1}(\l_{i,s+1}-\l^*_{i,s})c_i(x_{s+1})=k(-c(x_{s+1}),\l^*-\l_{s+1}).$$
 From concavity $d$ and $-c(x_{s+1})\in \partial d_y(\l_{s+1}) $ we obtain 
 \begin{equation}\label{Eq36_4.9.16}
 0\leq d_y(\l^*)-d_y(\l_{s+1})\leq (-c(x_{s+1}),\l^*-\l_{s+1}).
 \end{equation}
 From (\ref{Eq35_4.9.16}) and (\ref{Eq36_4.9.16}) follows 
 \begin{equation}\label{Eq37_4.9.16}
D(\l_s,\l^*)-D(\l_{s+1},\l^*)\geq k(d_y(\l^*)-d_y(\l_{s+1}))\geq 0.
 \end{equation}
 If $\lim_{s\to \infty}d_y(\l_s)=\bar{d}< d_y(\l^*)=F(y,x^*),$ then there is $\sigma >0$ and $s_0$ that from (\ref{Eq37_4.9.16}) we have
 $$D(\l_s,\l^*)-D(\l_{s+1},\l^*)\geq k\sigma, \forall s\geq s_0.$$
 
 Summing up the last inequalities from $s=s_0$ to $s=N$ we obtain
 $$D(\l_0,\l^*)-D(\l_{N+1},\l^*)=$$
 $$\sum^m_{i=1}(\l^*_i\ln\l_{i,N+1}/\l_{i,s_0}+\l_{i,s_0}-\l_{i,N+1})\geq k(N-s_0)\sigma,$$
 which is impossible for large $N$ due to the boundedness of $\{\l_s\}_{s=0}^{\infty}\in \Lambda_0$.
 
 Therefore 
 $$d_y(\l^*)=\lim d_y(\l_s)=\lim_{s\to\infty}[F(y,x_s)-(\l_s,c(x_s))].$$
 Keeping in mind asymptotic complementarity (\ref{Eq33_4.9.16}) we obtain
 \begin{equation}\label{Eq38_4.9.16}
 \lim_{s\to\infty}F(y,x_s)=F(y,x^*)=d_y(\l^*).
 \end{equation} 
 \item[3)]  From   boundedness of the dual sequence follows  existence of a subsequence $\{\l_{s_i}\}^{\infty}_{i=1}\subset\{\l_s\}^{\infty}_{s=0}$:$\lim_{s_i\to\infty}\l_{s_i}=\bar{\l}$. From    convergence of the dual sequence  in value follows $\bar{\l}=\l^*$ and $L^*=\{\l\in\Rm_+:d_y(\l)=d_y(\bar{\l})\}$.
 
 From dual strong monotonicity: $d_y(\l_{s+1})>d_y(\l_s)$ follows 
 $$L^*\subset... \Lambda_{s+1}\subset\Lambda_s...\subset \Lambda_0,$$
 therefore from (\ref{Eq26_4.9.16}) follows that $\{d_H(\partial \Lambda_s,L^*)\}^{\infty}_{s=0}$, is a monotone decreasing  sequence of positive numbers. It has a limit, that is
 $$\lim_{s\to\infty}d_H(\partial \Lambda_s,L^*)=\nu\geq 0,$$
 but $\nu>0$ is impossible due to the dual convergence in value (\ref{Eq38_4.9.16}).
 \item[4)] The ergodic convergence of the primal sequence one can prove by repeating the arguments used in the proof of item 4. Theorem 8 in \cite{polyak16}. The proof of Theorem \ref{T_29.11.6} is completed.

  \end{enumerate}
  So far, neither the fixed "center" $y\in \intr\;\Omega$ nor the fixed scaling parameter $k>0$ contributed to improvement of the EPM convergence.
   
   In the following section we establish $Q$-linear convergence rate of the EPM under standard second order sufficient optimality condition (\ref{Eq11_4.9.16})-(\ref{Eq12_4.9.16}), any  $y\in \intr\;\Omega$ as fixed "center" and a fixed, but large enough scaling parameter $k>0$. 
 %---------------------------------------------
 \section{Convergence rate of the EPM}

 Let us first describe  the dual domain, where the basic results are taking place.
 
 We assume that  $0<\delta<\min_{1\leq i\leq r}\l^*_i$ is small enough and $k_0>0$ is large enough.
 
 In the course of proving the Theorem it will be more clear quantitatively  what "small" and "large" means.
 
  We split the extended dual set into active and passive sub-sets, that is 
 $$\Lambda(\cdot)\equiv \Lambda(\l,k,\delta)=\Lambda_{(r)}(\cdot) \otimes \Lambda_{(m-r)}(\cdot),$$
 where 
 $$\Lambda_{(r)}(\cdot)\equiv \Lambda_{(r)}(\l_{(r)},k,\delta)=\{(\l_{(r)},k,\delta):\l_i\ge \delta,\; |\l_i-\l_i^*|\leq \delta k,i=1,...,r,k\geq k_0\}$$
 be the active dual sub-set and 
 $$\Lambda_{(m-r)}(\cdot)\equiv \Lambda_{(m-r)}(\l_{(m-r)}, k,\delta)$$ $$=\{(\l_{(m-r)}, k,\delta):0\leq \l_i\leq \delta k,\; i=r+1,...,m,\; k\geq k_0\}$$
be the passive dual sub-set.
 For a vector $a\in\Rn$ we use the following norm $\|a\|=\max_{1\leq i\leq n}|a_i|$. For a matrix $A\in \RR^{m\times n}$ the correspondent  norm is $\|A\|=\max_{1\leq i\leq m}\sum^n_{j=1}|a_{ij}|$.
 \begin{theorem}\label{T6_4.11.16}
 If $f$, $c_i\in C^2$, $i=1,...,m$ and the second order  sufficient optimality condition (\ref{Eq11_4.9.16})-(\ref{Eq12_4.9.16}) is satisfied, then exist  a small enough $\delta>0$ and large enough $k_0>0$, that for any  $y\in\intr\;\Omega$ as a fixed "center" and any $(\l,k)\in \Lambda(\cdot)$ the following statements hold true:
 \begin{enumerate}
\item[1)] there exists 
 $$\hat{x}=\hat{x}(\l,k):=\nabla_x\mathcal{L}_y(\hat{x},\l,k)=0$$
 and
 $$\hat{\l}=(\hat{\l}_i=\l_i(kc_i(\hat{x})+1)^{-1},\; i=1,...,m).$$
 \item[2)] for $(\hat{x},\hat{\l})$ the following bound holds
 \begin{equation}\label{Eq39_4.9.16}
 \max\{\|\hat{x}-x^*\|,\|\hat{\l}-\l^*\|\}\leq ck^{-1}\|\l-\l^*\|,
 \end{equation}
where  $c>0$ is independent on $k\geq k_0$. Also  $x(\l^*,k)=x^*$ and $\hat{\l}(\l^*,k)=\l^*$, that is $\l^*$ is a fixed point of the map $\l\to \hat{\l}(\l,k)$.
 \item[3)] The EDF $\mathcal{L}_y(x,\l,k)$ is strongly convex in the neighborhood of $\hat{x}$.
 \end{enumerate}
  \end{theorem}
 \proof Vector $\hat{x}=\hat{x}(\l,k)$ also  depends on $y \in \intr\;\Omega$, therefore $\hat{\l}$ is  a function of $\l,k$ and $y\in \intr\;\Omega$. At this point $y$ is fixed, so to simplify notation we omit $y$ in the definition of $\hat{x}$ and $\hat{\l}$.  
 
 By introducing vector $t=(t_1,...,t_r,t_{r+1},...,t_m)^T$ with $t_i=(\l_i-\l_i^*)k^{-1}$ we transform the  dual set $\Lambda(\cdot)$ into the following  neighborhood of the origin of the extended dual space
 $$S(0,k,\delta)=S_{(r)}(0,k,\delta)\otimes S_{(m-r)}(0,k,\delta),$$
 where
 $$S_{(r)}(0,k,\delta)=\{(t_{(r)},k):|t_i|\leq \delta,\; t_i\geq (\delta-\l_i^*)k^{-1},\; i=1,...,r,\ k\geq k_0\}$$
 and
 $$S_{(m-r)}(0,k,\delta)=\{(t_{(m-r)},k):0\leq t_i\leq \delta, i=r+1,...,m, k\geq k_0\}.$$
 Let us consider  vector-function  $h:\RR^{n+m-r+1}\to\Rn$ defined as follows
 $$h(x,t_{(m-r)},k)=k\sum^m_{i=r+1}t_i(kc_i(x)+1)^{-1}\nabla c_i(x),$$
 then
 $$\nabla_t h(x,t_{(m-r)},k)=[0^{n,r} \quad k\nabla c^T_{(m-r)}(x)\Psi^{'}(kc_{(m-r)}(x)],$$
 where $\Psi^{'}(kc_{(m-r)}(x))=\diag[(kc_i(x)+1)^{-1}]^m_{i=r+1}$.
 $$\nabla_xh(x,t_{(m-r)},k)=k^2\sum^m_{i=r+1}t_i(kc_i(x)+1)^{-2}\nabla c_i^T(x)\nabla c_i(x).$$
 Therefore for any $k>0$ we have 
 $$h(x^*,0^{m-r},k)=0^n,\; \nabla_x h(x^*,0^{m-r},k)=0^{n\times n}.$$
 Our main tool is the map $\Phi_y:\RR^{n+m+r+1}\to\RR^{n+r}$ given by the following formula
 $$\Phi_y(x,\hat{\l}_{(r)},t,k)=\begin{pmatrix}
  \Delta^{-1}(y,x)\nabla f(x)-\sum^r_{i=1}\hat{\l}_i\nabla c_i(x)-h(x,t_{(m-r)},k) \\
  (t_i+k^{-1}\l_i^*)(kc_i(x)+1)^{-1}-k^{-1}\hat{\l}_i,\; i=1,...,r  \\
  \end{pmatrix} .$$
  For a given fixed $y \in\intr\;\Omega$  we have
  $$\Phi_y (x^*,\l^*_{(r)},0^m,k)=\left[     \begin{array}{c}
        \Delta^{-1}(y,x^*)\nabla f(x^*)-\sum^r_{i=1}\l_i^*\nabla c_i(x^*)-h(x^*,0^{m-r},k) \\
        k^{-1}(\l^*_i-\l^*_i),\; i=1,...,r
      \end{array}
    \right]$$
  $$=\left[     \begin{array}{c}
         \nabla_x L_y(x^*,\l^*) \\
          k^{-1}(\l^*_i-\l^*_i),\;i=1,...,r
        \end{array}
      \right]=\left[     \begin{array}{c}
              0^n \\
              0^r
            \end{array}
          \right].$$
          Let us consider the following Jacobian
 $$\nabla_{x\hat{\l}_{(r)}}\Phi_y(x,\hat{\l}_{(r)},t,k)=\nabla_{x\hat{\l}_{(r)}}\Phi_y(\cdot)=$$
  $$= \left[     \begin{array}{cc}
               \nabla^2_{xx}L_y(\cdot) &- \nabla c^T_{(r)}(\cdot) \\
                -(T^r+k\Lambda^*_{(r)})(kC_{(r)}(\cdot)+I^r)^{-2}\nabla c_{(r)}(\cdot) &-k^{-1}I^r
              \end{array}
            \right],$$
where $\Lambda_{(r)}=\diag(\l_i)^r_{i=1}$, $C_{(r)}(\cdot)=\diag( c_i(\cdot))^r_{i=1}$, $T^r=\diag (t_i)_{i=1}^r$
$I^r$ - identical matrix in $\RR^r$.
For $x=x^*$, $\l_{(r)}=\l^*_{(r)}$ and $t=0^m$ we have
 $$\nabla_{x\hat{\l}_{(r)}}\Phi_y(x^*,\l^*_{(r)},0^m,k)= \left[     \begin{array}{cc}
                \nabla^2_{xx}L_y(x^*,\l^*) &- \nabla c^T_{(r)}(x^*) \\
                 -\Lambda^*_{(r)}\nabla c_{(r)}(x^*) &-k^{-1}I^r
               \end{array}
             \right] \equiv$$
             
             $$\equiv \left[     \begin{array}{cc}
                             \nabla^2_{xx}L_y &- \nabla c^T_{(r)} \\
                              -\Lambda^*_{(r)}\nabla c_{(r)} &-k^{-1}I^r
                            \end{array}
                          \right] \equiv\nabla \Phi_{(y,k)}.$$

 The next step is to show that the matrix $\nabla\Phi_{(y,k)}$ is not singular for any given $y\in \intr\;\Omega$ and $k\geq k_0$, where $k_0>0$ is large enough.
 
 Let $w=(u,v)\in\RR^{n+r}$, then from
 $$\nabla \Phi_{(y,k)}w=  \left[     \begin{array}{cc}
                 \nabla^2_{xx}L_yu &- \nabla c^T_{(r)}v \\
                  -\Lambda^*_{(r)}\nabla c_{(r)}u &-k^{-1}v
                \end{array}
              \right]= \left[     \begin{array}{c}
                            0^n \\
                            0^r
                          \end{array}
                        \right] ,$$
follows $v=-k\Lambda^*_{(r)}\nabla c_{(r)}u $ and 
$$Nu=(\nabla^2_{xx}L_y+k\nabla c^T_{(r)}\Lambda^*_{(r)}\nabla c_{(r)})u=0,$$
  therefore
$$(Nu,u)=((\nabla^2_{xx}L_y+k\nabla c_{(r)}^T\Lambda^*_{(r)}\nabla c_{(r)})u,u)=0.$$
  From sufficient optimality condition (\ref{Eq11_4.9.16})-(\ref{Eq12_4.9.16}) and Debreu's lemma with $A=\nabla^2_{xx}L_y$ and $C=\Lambda^{*\frac{1}{2}}\nabla c_{(r)}$ follows the existence $0<\rho<\mu$ that
  $$0=(Nu,u)\geq\rho(u,u),$$
  hence $u=0^n$, then from 
  $$\nabla^2_{xx}L_yu-\nabla c^T_{(r)}v=0^n$$
  and (\ref{Eq12_4.9.16}) follows $v=0^r$. It means 
  $$\nabla \Phi_{(y,k)}w=0^{n+r}\Rightarrow w=0^{n+r},$$
  therefore $\nabla \Phi^{-1}_{(y,k)}$ exists. Using argument similar to those in Theorem 1 from  \cite{polyak92} one can prove existence  of  large enough $k_0>0$, that for any $k\geq k_0$ there exists $\rho_0>0$ independent on $k\geq k_0$ and  $y\in \intr\;\Omega$ that 
  \begin{equation}\label{Eq40_4.9.16}
  \|\Phi^{-1}_{(y,k)}\|\leq \rho_0.
  \end{equation}
  Let $\infty>k_1>k_0$, $k_0>0$ be large enough and $K=\{0^n\}\times[k_0,k_1]$.  We consider the following neighborhood  
  $$S(K,\delta)=\{(t,k):|t_i|\leq\delta,\; t_i\geq (\delta-\l^*_i)k^{-1},i=1,...,r;$$ $$0\leq t_i\leq \delta, i=r+1,...,m\}$$
  of $K$. 
  
  From the second implicit function Theorem  (see, for example,  \cite{bertsekas82} p.12) follows that for any $k\in [k_0,k_1]$ the system
  $$\Phi_y (x,\hat{\l}_{(r)},t,k)=0^{m+r}$$
  defines on $S(K,\delta)$  a unique pair of vectors $$x(t,k)=(x_i(t,k),\;i=1,...,n) \mbox{ and  } \hat{\l}_{(r)}(t,k)=(\hat{\l}_i(t,k),i=1,...,r)$$   that $ x(0^m,k)=x^*$, $\hat{\l}_{(r)}(0^m,k)=\l^*_{(r)}$ and
  \begin{equation}\label{Eq41_4.9.16}
  \Phi_y(x(t,k),\hat{\l}_{(r)}(t,k),t,k)\equiv0^{n+r},\forall (t,k)\in S(K,\delta).
  \end{equation} 
  
  Identity (\ref{Eq41_4.9.16}) can be rewritten as follows
  \begin{equation}\label{Eq42_4.9.16}
  \Delta^{-1}(y,x(t,k))\nabla f(x(t,k))-\sum^r_{i=1}\hat{\l}_i(t,k)\nabla c_i(x(t,k))
  \end{equation}
  $$-h(x(t,k),t_{(m-r)},k)\equiv 0^n.$$
 \begin{equation}\label{Eq43_4.9.16}
 \hat{\l}_i(t,k)\equiv (kt_i+\l^*_i)\psi^{'}(kc_i(x(t,k))),\; i=1,...,r,
 \end{equation}
 also 
 \begin{equation}\label{Eq44_4.9.16}
 \hat{\l}_i(t,k)\equiv kt_i\psi^{'}(kc_i(x(t,k))),\; i=r+1,...,m.
 \end{equation}
 From (\ref{Eq42_4.9.16})-(\ref{Eq44_4.9.16}) follows 
 $$\nabla_x\mathcal{L}_y(x(t,k),\l,k)=\nabla_xL_y(x(t,k),\hat{\l}(t,k))\equiv 0^n,$$
 where 
 $\hat{\l}(t,k)=(\hat{\l}_{(r)}(t,k),\hat{\l}_{(m-r)}(t,k))^T.$
 
 It completes the proof of item 1).

 2)  For a given small enough  $\varepsilon>0$ there is $\delta>0$ such that 
 $$\|x(t,k)-x(0^m,k)\|=\|x(t,k)-x^*\|\leq \varepsilon \mbox{ for }  \forall (t,k)\in S(K,\delta).$$
  Hence, there is $\sigma=\min_{r+1 \leq i\leq m}\{c_i(x^*)\}$ that  for $r+1 \leq i\leq m$ we have
  $$c_i(x(t,k))\geq 0.5 \sigma,\;\;\; \forall (t,k) \in S(K,\delta).$$
  From (\ref{Eq17_4.9.16}) follows
  $$\hat{\l}_i=\l_i(kc(x(t,k)+1))^{-1}\leq 2(\sigma k)^{-1}\l_i,\; i=r+1,..., m,$$
  where  $\sigma>0$ is independent on $k\geq k_0$. 
  
  To prove the bound (\ref{Eq39_4.9.16}) for $x(t,k)$ and $\hat{\l}_{(r)}(t,k)$ we will first estimate the norms $\|\nabla_t x(t,k)\|$, $\|\nabla_t\nabla{\hat{\l}}_{(r)}(t,k)\|$   at  $t=0^m$.

  By differentiating identities (\ref{Eq42_4.9.16}) and (\ref{Eq43_4.9.16}) in $t$ we obtain the following  system for 
   Jacobians $\nabla_tx(\cdot)$ and $\nabla_t\hat{\l}_{(r)}(\cdot)$
  \begin{equation}\label{Eq45_4.9.16}
  \nabla^2_{xx}L_y(x(\cdot),\hat{\l}_{(r)}(\cdot))\nabla_tx(\cdot)-\nabla c_{(r)}^T(\cdot)\nabla_t\hat{\l}_{(r)}(\cdot)
  \end{equation}
  $$\equiv \nabla_t h(x(t,k),t_{(m-r)},k)$$
  $$=\nabla_xh(x(\cdot),\cdot)\nabla_t(x(\cdot))+\nabla_th(x(\cdot),\cdot)$$
 \begin{equation}\label{Eq46_4.9.16}
 k(kT^r+\Lambda^*_{(r)})\Psi^{''}(kc_{(r)}(x(\cdot)))\nabla c_{(r)}(x(\cdot))\nabla_t(x(\cdot))-\nabla_t\hat{\l}_{(r)}(\cdot)
 \end{equation}
 $$\equiv -[k\Psi^{'}(kc_{(r)}(x(\cdot))); 0^{r,m-r}],$$
 
  where $\Psi^{'}(kc_{(r)}(x(\cdot)))=\diag (\psi^{'}(kc_i(x(\cdot))))^r_{i=1}$, $\Psi^{''}(kc_{(r)}(x(\cdot)))=\diag(\psi^{''}(kc_i(x(\cdot))))^r_{i=1}$  are diagonal matrices. 
  
 In other words, for Jacobians $\nabla_tx(\cdot)$ and $\nabla_t\hat{\l}_{(r)}(\cdot)$ we have the following system
 \begin{equation}\label{Eq47_4.9.16}
 \nabla\Phi_{(y,k)}(\cdot)\left[\begin{array}{c}
 \nabla_tx(\cdot)\\ \nabla_t\hat{\l}(\cdot)
 \end{array}\right]=R(\cdot),
 \end{equation}
 where 
 $$\nabla\Phi_{(y,k)}(\cdot)=\left[\begin{array}{cc }
 \nabla^2_{xx}L_y(x(\cdot),\hat{\l}_{(r)}(\cdot))& -\nabla c^T_{(r)}(x(\cdot))\\
  (kT^r+\Lambda^{*}_{(r)})\Psi^{''}(kc_{(r)}(x(\cdot)))\nabla c_{(r)}(x(\cdot))&-k^{-1}I^r
 \end{array}\right]$$
 and
 $$R(x(t,k);t,k)=R(x(\cdot),\cdot)=\left[\begin{array}{c}
 \nabla_xh(x(\cdot),\cdot)\nabla_tx(\cdot)+\nabla_th(x(\cdot),\cdot)\\\Psi^{'}(kc_{(r)}(x(\cdot))); 0^{r,m-r}
 \end{array}\right].$$
Let us consider the system (\ref{Eq47_4.9.16})  for $t=0^m$. We obtain
$$x(0^m,k)=x^*,\; \hat{\l}_{(r)}(0^m,k)=\l^*_{(r)},$$
$$\nabla^2_{xx}L_y(x(0^m,k),\hat{\l}_{(r)}(0^m,k))=\nabla^2_{xx}L_y(x^*,\l^*)=\nabla^2_{xx}L_y,$$
$$\nabla c_r(x(0^m,k))=\nabla c_{(r)}(x^*)=\nabla c_{(r)}$$
$$\Psi^{'}(kc_{(r)}(x(0^m,k)))=\Psi^{'}(kc_{(r)}(x^*))=I^r,\; \Psi^{''}(kc_{(r)}(x^*))=\psi^{''}(0)I^r.$$
We also have 
$$\nabla_th(x(t,k),t,k)_{|t=0^m}\equiv \nabla_th(x(\cdot),\cdot)_{|t=0^m}$$
$$=[0^{n,r}\;\; k\nabla c^T_{(m-r)}(x^*)\Psi^{'}(kc_{(m-r)}(x^*))],$$
where $\Psi^{'}(kc_{(m-r)}(x^*))=\diag[\psi^{'}(kc_i(x^*))]^m_{i=r+1}$.

From $c_i(x^*)\geq \sigma>0$, $r+1\leq i \leq m$ and the update formulas (\ref{Eq17_4.9.16}) follows
$$\|k(\nabla c_{(m-r)}(x^*))^T\Psi^{'}(kc_{(m-r)}(x^*))\|\leq 2\sigma ^{-1}\|(\nabla c_{(m-r)}(x^*))^T\|.$$
The system (\ref{Eq47_4.9.16}) for $t=0^m$ we can rewrite as follows
$$\left[     \begin{array}{cc}
                            \nabla^2_{xx}L_y &-\nabla c^T_{(r)} \\
                            -\Lambda^*_{(r)}\nabla c_{(r)} &-k^{-1}I^r
                          \end{array}
                        \right]\left[     \begin{array}{c}
                                                    \nabla_t x(0^m,k) \\
                                                   \nabla_t\hat{\l}_{(r)}(0^m,k)
                                                  \end{array}
                                                \right]=\left[     \begin{array}{cc}
                                                                            0^{n,r} & k(\nabla c_{(m-r)}(x^*))^T\Psi^{'}(kc_{(m-r)}(x^*)\\
                                                                            I^r & 0^{r,m-r}
                                                                          \end{array}
                                                                        \right] .$$
   or
   \begin{equation}\label{Eq48_4.9.16}
   \left[\begin{array}{c}
                               \nabla_t x(0^m,k) \\
                               \nabla_t\hat{\l}_{(r)}(0^m,k)
                             \end{array}
                           \right] 
   \end{equation}        
   $$=\nabla\Phi^{-1}_{(y,k)}\left[ \begin{array}{cc}
                                                          0^{n,r} & k(\nabla_{(m-r)} c(x^*))^T\Psi^{'}(kc_{(m-r)}(x^*))\\
                                                          I^r & 0^{r,m-r}
                                                        \end{array}
                                                      \right].$$                                                             
From (\ref{Eq40_4.9.16}), \eqref{Eq48_4.9.16} and $k\geq k_0$ follows 
$$\max\{\|\nabla_t\hat{x}(0^m,k)\|,\|\nabla_t\hat{\l}_{(r)}(0^m,k)\|\}\leq \rho_0 \max\{1, 2\sigma^{-1}\|(\nabla c_{(m-r)}(x^*))^T\|\}=c_0.$$
Thus, for $\delta>0$ small enough and any $(t,k)\in S(K,\delta)$ from \eqref{Eq47_4.9.16} follows 
\begin{equation}\label{Eq49_4.9.16}
\|\Phi^{-1}_{(y,k)}(\cdot)R(\cdot)\|=
\end{equation}
$$=\|\nabla\Phi^{-1}_{(y,k)}(x(\tau t,k),\hat{\l}_{(r)} (\tau t, k))\cdot R (x(\tau t, k); \tau t, k)\|\leq 2 c_0$$
for any $0\leq \tau\leq 1$ and any $k\geq k_0$.

Using Newton-Leibniz formula
\begin{equation}\left[     \begin{array}{c}
                            x(t,k)-x^* \\
                            \hat{\l}_{(r)}(t,k)-\l^*_{(r)}
                          \end{array}
                        \right] =\left[     \begin{array}{c}
                                                    x(t,k)-x(0^m,k) \\
                                                    \hat{\l}_{(r)}(t,k)-\hat{\l}_{(r)}(0^m,k)
                                                  \end{array}
                                                \right] \end{equation}
 $$\label{Eq50_4.9.16}=\int^1_0\nabla \Phi^{-1}_{(y,k)}\left((x(\tau t, k), \hat{\l}_{(r)}(\tau t, k)\right)R(x(\tau t, k),\tau t,k)[t]d\tau
 $$
 we obtain  
 $$\max\{\|x(t,k)-x^*\|,\|\hat{\l}_{(r)}(t,k)-\l^*\|\}\leq 2c_0\| t\|=2c_0k^{-1}\|\l-\l^*\|.$$
 Let  $\hat{x}(\l,k)=x(\frac{\l-\l^*}{k},k)$ and 
 $$\hat{\l}(\l,k)=(\hat{\l}_{(r)}(\frac{\l-\l^*}{k}),k), \hat{\l}_{(m-r)}(\frac{\l-\l^*}{k},k)).$$
 Then for $c=2\max\{\sigma^{-1}, c_0\}$, which is independent on $k\geq k_0$, we obtain \eqref{Eq39_4.9.16} 
 
 3) Let us consider the Hessian of the EDF $\mathcal{L}_y(x,\l,k)$ at $x=\hat{x}$ and $\l=\hat{\l}$. 
 We have 
 $$\nabla^2_{xx}\mathcal{L}_y(\hat{x},\hat{\l},k)=\nabla^2_{xx}L_y(\hat{x},\hat{\l})-k\nabla c(\hat{x})^T\Phi^{''}(kc(\hat{x}))\hat{\Lambda}\nabla c(\hat{x}).$$
 From (\ref{Eq39_4.9.16}) for $k\geq k_0$ large enough we have 
 $$\nabla^2_{xx}\mathcal{L}_y(\hat{x},\hat{\l},k)\approx \nabla^2_{xx}L_y(x^*,\l^*)-k\psi^{''}(0)\nabla c(x^*)^T\Lambda^*\nabla c  (x^*)$$
 $$\nabla^2_{xx}L_y(x^*,\l^*)+k\nabla c_{(r)}(x^*)^T\Lambda^*_{(r)}\nabla c_{(r)}(x^*).$$
 The  item 3) of Theorem \ref{T6_4.11.16} follows from the second order sufficient optimality condition (\ref{Eq11_4.9.16})-(\ref{Eq12_4.9.16}) and Debreu's  Lemma, with $A=\nabla^2_{xx}L_y(x^*,\l^*)$ and $C=\Lambda^{*\frac{1}{2}}_{(r)}\nabla c_{(r)}(x^*)$.
 The proof of Theorem \ref{T6_4.11.16} is completed.
 %----------------------------------------------------------------------------------
 \section{Stopping Criteria}
The EPM (\ref{Eq16_4.9.16})-(\ref{Eq17_4.9.16}) is an infinite procedure, which require, at each step, solving an unconstrained optimization problem (\ref{Eq16_4.9.16}), which is, in turn,  an infinite procedure as well.
 
 The following result allows replacing  $x_{s+1}$ from (\ref{Eq16_4.9.16}) by an approximation $\bar{x}_{s+1}$, finding  which requires finite procedure and does not compromise $Q$-linear convergence rate.
 
 For a given $\alpha>0$ let us consider the primal-dual approximation $(\bar{x},\bar{\l}):
$
\begin{equation}\label{Eq_4.11.16}
\bar{x}=\bar{x}(\l,k):\|\nabla_x\mathcal{L}_y(\bar{x},\l,k)\|\leq\frac{\alpha}{k}\|\bar{\l}-\l\|
\end{equation}
\begin{equation}\label{Eq2_4.11.16}
\bar{\l}=\bar{\l}(\l,k)=(\bar{\l}_i=\psi^{'}(kc_i(\bar{x}))\l_i,\;i=1,...,m.)
\end{equation}
 Obviously $\bar{x}$ depends not only on $\l\in\Rm_{++}$ and $k>0$ but also on $y$ and $\alpha$ as well.
 
 At this point $y\in \intr\;\Omega$ and $\alpha>0$ are fixed, therefore to simplify notation we omitted $y$ and $\alpha$ from the definition of $\bar{x}$ and $\bar{\l}$.
 
 \begin{theorem}\label{T_27.11.16}
 If $f$, $c_i\in C^2$,   $i=1,...,m$ and the second order sufficient optimality condition (\ref{Eq11_4.9.16})-(\ref{Eq12_4.9.16}) is satisfied, then for a given $\alpha>0$, small enough $\delta>0$, large enough $k_0$, any $k\geq k_0$ and any $(\l,k)\in\Lambda(\l,k,\delta)$ we have:
 \begin{enumerate}
\item[1)]  there exists $(\bar{x},\bar{\l})$ defined by (\ref{Eq_4.11.16})-(\ref{Eq2_4.11.16});
\item[2)] there is $c>0$  independent on $k\geq k_0$ that the following bound
\begin{equation}\label{Eq3_4.11.16}
\max\{\|\bar{x}-x^*\|,\|\bar{\l}-\l^*\|\}\leq \frac{c}{k}(1+2\alpha)\|\l-\l^*\|
\end{equation}
holds;
\item[3)] the Lagrangian $\mathcal{L}_y(x,\l,k)$  for the equivalent problem is strongly convex at the neighborhood of $\bar{x}$.
 \end{enumerate}
 \end{theorem}
 \proof For a small enough $\delta>0$, large enough $k_0$ and any $k\geq k_0$ we define the following extended dual set 
 $$\Lambda(\l,k,\delta,\theta)=\Lambda(\l,k,\delta)\otimes\{\theta\in\Rn:\|\theta\|\leq \delta\}$$
 $$=\{\l\in\Rm_+:\l_i\geq\delta,|\l_i-\l^*_i|\leq\delta k,\;i=1,...,r\}\otimes\{0<\l_i<\delta k, \;i=r+1,...,m,\;k\geq k_0\}$$
 $$\otimes\{\theta\in\Rn:\|\theta\|\leq \delta\}=\Lambda_{(r)}(\cdot)\otimes\Lambda_{(m-r)}(\cdot)\otimes \{\theta\in\Rn:\|\theta\|\leq \delta\}.$$
 
 By introducing vector $t=(t_1,...,t_r,t_{r+1},...,t_m)$ with $t_i=(\l_i-\l^*_i)k^{-1}$, $i=1,...,m$ we transform $\Lambda(\l,k,\delta, \theta)$ into the neighborhood of the origin of the extended dual space
 $$S(0^m,k,\delta,0^n)=S_{(r)}(0^r,k,\delta)\otimes S_{(m-r)}(0^{m-r},k,\delta)\otimes\{\theta\in\Rn:\|\theta\|\leq\delta\}.$$
 The following map $\Phi_y:\RR^{2n+m+r+1}\to\RR^{n+r}:$
 $$\Phi_y(x,\bar{\l}_{(r)},k,t,\theta) =\left(     \begin{array}{c}
                                                   \Delta^{-1}(y,x)\nabla f(x)-\sum^r_{i=1}\bar{\l_i}\nabla c_i(x)-h(x,t_{(m-r)},k)-\theta \\
                                                    (t_i+k^{-1}\l_i^*)\psi^{'}(kc_i(x))-k^{-1}\bar{\l}_i,\;\;i=1,...,r
                                                   \end{array}
                                                 \right) ,$$
 is the key ingredient of the proof.

 Let us consider $\infty>k_1>k_0$, then for any $k\in[k_0,k_1]$ and $y\in \intr\;\Omega$ we have 
$$ \Phi_y(x^*,\l^*_{(r)},k,0^m,0^n)=$$ $$\left[     \begin{array}{c}
                           \Delta^{-1}(y,x^*)\nabla f(x^*)-\sum^r_{i=1}\l^*_i\nabla c_i(x^*)-h(x^*,0^{m-r},k)-0^n \\
                            k^{-1}(\l^*_i-\l^*_i),\;\; i=1,...,r
                           \end{array}
                         \right] =\left[     \begin{array}{c}
                                                    0^n \\
                                                    0^r
                                                   \end{array}
                                                 \right] .$$
 Further,
$$\nabla_{x\bar{\l}_{(r)}}\Phi_{y}\equiv \nabla_{x,\bar{\l}_{(r)}}\Phi_y(x^*,\l^*_{(r)},k,0^m,0^n)=$$
$$ \left[     \begin{array}{c}
                            \nabla^2_{xx}L_y(x^*,\l^*)-\nabla c^T_{(r)}(x^*) \\
                             \Lambda^*_{(r)}\nabla c_{(r)}(x^*)-k^{-1}I^r
                           \end{array}
                         \right]=\nabla \Phi_{(y,k)}.$$
 
We saw already that the inverse 
$(\nabla \Phi_{(y,k)})^{-1}$ exists and there is $\rho_0>0$ that $\|(\nabla\Phi_{(y,k)})^{-1}\|\leq \rho_0$.

From the second Implicit Function Theorem follows existence of two vector-functions 
$$x(\cdot)=x(t,k,\theta)=(x_1(t,k,\theta),...,x_n(t,k,\theta))$$
and
$$\bar{\l}_{(r)}(\cdot)=\bar{\l}_{(r)}(t,k,\theta)=(\bar{\l}_1(t,k,\theta),...,\bar{\l}_r(t,k,\theta))$$
uniquely defined on $S(0^m,k,\delta, 0^n)$ for small enough $\delta>0$ and $k\geq k_0$, that the following identities hold
\begin{equation}\label{Eq4_4.11.16}
\Delta^{-1}(y,x(\cdot))\nabla f(x(\cdot))-\sum^r_{i=1}\bar{\l}_i\nabla c_i(x(\cdot))-h(x(\cdot),t_{(m-r)},k)-\theta\equiv 0
\end{equation}
\begin{equation}\label{Eq5_4.11.16}
\bar{\l}_i(\cdot)=\bar{\l}_i(t,k,\theta)\equiv (kt_i+\l^*_i)\psi^{'}(kc_i(x(\cdot))),\;\;i=1,..,r.
\end{equation}
For a given small $\varepsilon>0$ there is $\delta>0$ that 
$$\max\{\|x(t,k,\theta)-x^*\|,\|\bar{\l}_{(r)}(t,k,\theta)-\l^*_{(r)}\|\}\leq \varepsilon$$
for $\forall (t,k,\theta)\in S(0^m,k,\delta,0^n)$.

Therefore for the passive constraints we have 
$$c_i(x(\cdot))=\bar{c}_i(x(t,k,\theta))\geq 0.5\sigma.$$
Hence
$$\bar{\l}_i=\l_i\psi^{'}(kc_i(x(\cdot)))\leq\l_i\psi^{'}(0.5k\sigma)\leq \frac{2}{\sigma k}\l_i,\;i=r+1,...,m$$
and $\sigma>0$ is  independent on $k\in [k_0,k_1]$.

To prove the bound (\ref{Eq3_4.11.16}) we estimate the norms of  Jacobians $\nabla_{t,\theta} x(t,k,\theta)$ and $\nabla_{t,\theta}\bar{\l}_{(r)}(t,k,\theta)$ for $t=0^m$, $\theta=0^n$.

By differentiating identities (\ref{Eq4_4.11.16}) and (\ref{Eq5_4.11.16}) in $t$ and $\theta$ we obtain
\begin{equation}\label{Eq6_4.11.16}
\nabla^2_{xx}L_y(\cdot)\nabla_{t,\theta}x(\cdot)-\nabla c^T_{(r)}(x(\cdot))\nabla_{t,\theta}\bar{\l}_{(r)}(\cdot)
\end{equation}
$$=[0^{n,r},\;\;\nabla_x h(x(\cdot),t_{(m-r)},k)\nabla_{t,\theta}x(\cdot)+\nabla_{t}h(x(\cdot),t_{(m-r)},k), \;\; I^n]$$
\begin{equation}\label{Eq7_4.11.16}
k(kT^r+\Lambda^*_{(r)})\Psi^{''}(kc_{(r)}(x(\cdot)))\nabla c_{(r)}(x(\cdot))\nabla_{t,\theta}x(\cdot)-\nabla_{t,\theta}\bar{\l}_{(r)}(\cdot)
\end{equation}
$$=-[k\Psi^{'}(kc_{(r)}(x(\cdot))), \;\;0^{r,m-r},\;\; 0^{r,n}],$$
where $T^r=\diag(t_i)^r_{i=1}$, $\Lambda^*_{(r)}=\diag(\l^*_i)^r_{I=1},$
$\Psi^{'}(kc_{(r)}(x(\cdot)))=\diag(\psi^{'}(kc_i(x(\cdot))))^r_{i=1}$, 
$\Psi^{''}(kc_{(r)}(x(\cdot))=\diag(\psi^{''}(kc_i(x(\cdot))))^r_{i=1}$.

The system (\ref{Eq6_4.11.16})-(\ref{Eq7_4.11.16}) can be rewritten as follows 
  \begin{equation}\label{Eq8_4.11.16}\left[     \begin{array}{c}
                          \nabla_{t,\theta}x(\cdot) \\
                           \nabla_{t,\theta}\bar{\l}_{(r)}(\cdot)
                           \end{array}
                         \right] =
                \end{equation}
                          $$\left[     \begin{array}{cc}
                                                   \nabla^2L_y(x(\cdot),\bar{\l}_{(r)}(\cdot)) & -\nabla c^T_{(r)}(\cdot) \\
                                                    (kT^r+\Lambda^*_{(r)})\Psi^{''}(kc_{(r)}(x(\cdot))\nabla c_{(r)}(x(\cdot))) & -k^{-1}I^r
                                                   \end{array}
                                                 \right]^{-1}\times  R(x(\cdot),t,k) ,$$
where 
$$ R(x(\cdot),t,k)=\left[     \begin{array}{ccc}
                                                   0^{n,r} & \nabla_xh(x(\cdot),t_{(m-r)},k)\nabla_{t,\theta}x(\cdot)+\nabla_th(x(\cdot),t_{(m-r)},k) & I^n \\
                                                    \Psi^{'}(kc_{(r)}(x(\cdot))) & 0^{r,m-r} & 0^{r,n}
                                                   \end{array}
                                                 \right].
                                                $$
    Let us consider system (\ref{Eq8_4.11.16}) for $t=0^m$ and $\theta=0^n$. We obtain 
\begin{equation}\label{Eq9_4.11.16}
 \left[     \begin{array}{c}
                            \nabla_{t,\theta}x(0^m,k,0^n) \\
                            \nabla_{t,\theta}\bar{\l}_{(r)}(0^m,k,0^n)
                           \end{array}
                         \right]= \left[     \begin{array}{cc}
                                                     \nabla^2 L_y(x^*,\l^*_{(r)}) & -\nabla c^T_{(r)}(x^*) \\
                                                      -\Lambda^*_{(r)}\nabla c_{(r)}(x^*) & -k^{-1}I^r
                                                    \end{array}
                                                  \right]^{-1}\times 
\end{equation}
$$ \left[     \begin{array}{ccc}
                         0^{n,r} & k\nabla c^T_{(m-r)}(x^*)\Psi^{'}(kc_{(m-r)}(x^*) & I^n \\
                            I^r & 0^{r,m-r} & 0^{r,n}
                           \end{array}
                         \right]=\Phi^{-1}_{(y,k)}R,$$
                         where $$\Psi^{'} (kc_{(m-r)}(x^*))=\diag (\psi^{'}(kc_i(x^*)))^m_{i=r+1}.$$
                         
Keeping in mind (\ref{Eq40_4.9.16}) from (\ref{Eq9_4.11.16}) we obtain 
$$\max\{\|\nabla_{t,\theta}x(0^m,k,0^n)\|,\; \|\nabla_{t,\theta}\bar{\l}_{(r)}(0^m,k,0^n)\|\}\leq$$
$$ \rho_0\max\{1,\sigma^{-1}\|(\nabla c_{(m-r)}(x^*))^T\|\}=c_0$$ and $c_0$  is independent on $k\geq k_0$.
Thus, 
$$ \left[     \begin{array}{c}
                            x(t,k,\theta)-x^*\\
                             \bar{\l}_{(r)}(t,k,\theta )-\l^*_{(r)}
                           \end{array}
                         \right]=  \left[     \begin{array}{c}
                                                     x(t,k,\theta)-x(0^m,k,0^n) \\
                                                    \bar{\l}_{(r)}(t,k,\theta)-\bar{\l}_{(r)}(0^m,k,0^n)
                                                    \end{array}
                                                  \right]=$$
                                                  $$=\int^1_0\nabla\Phi^{-1}_{(y,k)}(x(\tau t,k,\tau\theta),\bar{\l}_{(r)}(\tau t,k, \tau \theta))R(x(\tau t,k,\tau\theta);\tau t,k,\tau\theta) \left[     \begin{array}{c}
                                                                              t \\
                                                                              \theta
                                                                             \end{array}
                                                                           \right]d\tau.$$
Hence, $$\max\{\|x(t,k,\theta)-x^*\|,\|\bar{\l}_{(r)}(t,k,\theta)-\l^*_{(r)}\|\}\leq 2c_0k^{-1}\|\l-\l^*\|+\|\theta\|.$$

Let 
$$\bar{x}=\bar{x}(t,k,\theta)=x\left(\frac{\l-\l^*}{k},k,\theta\right)$$
$$\bar{\l}=\bar{\l}(t,k,\theta)=\left(\bar{\l}_{(r)}\left(\frac{\l-\l^*}{k},k,\theta \right);\bar{\l}_{(m-r)}\left(\frac{\l-\l^*}{k},k,\theta\right)\right),$$

then for $c=2\max\{\sigma^{-1},c_0\}$ we have 
\begin{equation}
\|\bar{x}-x^*\|\leq\frac{c}{k}\|\l-\l^*\|+\|\theta\|
\end{equation}
\begin{equation}
\|\bar{\l}-\l^*\|\leq \frac{c}{k}\|\l-\l^*\|+\|\theta\|.
\end{equation}
Keeping in mind the stopping criteria  (\ref{Eq_4.11.16})-(\ref{Eq2_4.11.16}) we obtain
$$\|\nabla_x\mathcal{L}_y(\bar{x},\l,k)\|=\|\theta\|\leq\frac{\alpha}{k}\|\bar{\l}-\l\|.$$
Therefore
\begin{equation}\label{Eq12_4.11.16}
\|\bar{x}-x^*\|\leq\frac{c}{k}\|\l-\l^*\|+\frac{\alpha}{k}\|\bar{\l}-\l\|
\end{equation}
\begin{equation}\label{Eq13_4.11.16}
\|\bar{\l}-\l^*\|\leq\frac{c}{k}\|\l-\l^*\|+\frac{\alpha}{k}\|\bar{\l}-\l\|.
\end{equation}
From (\ref{Eq13_4.11.16}) follows  
$$\|\bar{\l}-\l^*\|\leq \frac{c}{k}\|\l-\l^*\|+\frac{\alpha}{k}\|\bar{\l}-\l^*\|+\frac{\alpha}{k}\|\l^*-\l\|$$
or
\begin{equation}\label{Eq14_4.11.16}
\left(1-\frac{\alpha}{k}\right)\|\bar{\l}-\l^*\|\leq\frac{c+\alpha}{k}\|\l-\l^*\|.
\end{equation}
 For  $k_0>c+2\alpha$ and any $k\geq k_0$ from (\ref{Eq14_4.11.16}) follows 
  \begin{equation}\label{Eq_18.11.16}
 \|\bar{\l}-\l^*\|\leq\frac{c+2\alpha}{k}\|\l-\l^*\|
 \end{equation}

From (\ref{Eq12_4.11.16}) and (\ref{Eq_18.11.16}) we  obtain
$$\|\bar{x}-x^*\|\leq\frac{c}{k}\|\l-\l^*\|+\frac{\alpha}{k}\|\bar{\l}-\l^*\|+\frac{\alpha}{k}\|\l-\l^*\|$$
$$=\frac{c+\alpha}{k}\|\l-\l^*\|+\frac{\alpha}{k}\|\bar{\l}-\l^*\|\leq\left[\frac{c+\alpha}{k}+\frac{\alpha}{k}\frac{(c+2\alpha)}{k}\right]\|\l-\l^*\|.$$
Again for $k_0>c+2\alpha$ and any $k\geq k_0$ we have 
$$\|\bar{x}-x^*\|\leq\frac{c+2\alpha}{k}\|\l-\l^*\|.$$
The proof of Theorem \ref{T_27.11.16} is completed.

We conclude the section by considering the numerical realization of the EPM.

The EPM scheme consists of inner and outer iteration.  On the inner iteration we find an approximation $\bar{x}$ for the primal  minimizer using the stopping criteria (\ref{Eq_4.11.16}).

On the outer iteration we update the Lagrange multipliers  by (\ref{Eq2_4.11.16}),  using the approximation $\bar{x}$.

For finding $\bar{x}$ any unconstrained minimization technique can be used. Fast gradient method ( see \cite{nesterov04}) or regularized Newton method (see \cite{polyak09}) are two possible candidates.

Under usual convexity and smoothness assumptions both methods converges to the minimizer from any starting point and for both methods there exist complexity bounds, that is the upper bound for the number of step required for finding an $\varepsilon$- approximation for the minimizer. 

To describe the numerical realization of EPM we need to introduce the relaxation operator $R:\Omega_{-k^{-1}}\times\RR^m_{++}\to \Omega_{-k^{-1}}\times \RR^m_{++}$, with is  defined as follows
\begin{equation}\label{Eq4_27.11.16}
Ru=\bar{u}=(\bar{x},\bar{\l}),
\end{equation}
where $\bar{x}$ and $\bar{\l}$ are given by \eqref{Eq_4.11.16} and \eqref{Eq2_4.11.16}.

We also need the merit function $\nu_y:\Omega_{-k^{-1}}\times\RR^m_{+}\to R$, which is defined by the following formula 
\begin{equation}\label{Eq_27.11.16}
\nu_y(u)=\max\{\|\nabla_xL_y(x,\l)\|,\;\;\sum^m_{i=1}\l_i|c_i(x)|, -c_i(x),\; i=1,...,m\}.
\end{equation}
From \eqref{Eq_27.11.16} follows $\nu_y(u)\geq 0$, $\forall u\in \Omega_{-k^{-1}}\times \RR^m_+$, it is also easy to see that 
\begin{equation}
\nu_y(u)=0\Leftrightarrow u=u^*=(x^*;\l^*)
\end{equation}
holds.

Moreover, under the second order sufficient optimality condition and $f,$ $c_i\in C^2$, $ i=1,...,m$  the merit function $\nu_y$ in the neighborhood of $u^*$ is similar to the norm of a gradient of a strong convex function with Lipschitz continuous gradient in the neighborhood of the minimizer (see \cite{polyak09_2}).

Let $\gamma>0$ be small enough, $y\in \intr\;\Omega_{\gamma}$ be the initial "center", $u=(x;\l)\in\Omega_{-k^{-1}}\times \RR^m_{++}$ be the initial primal-dual approximation, $\Delta>0$ be the reduction parameter for the objective function, $k>0$ be the scaling parameter and $\varepsilon>0$ be the required accuracy.

The EPM consists of the following operations
\begin{enumerate}
\item[1.] find $\bar{u}=(\bar{x};\bar{\l})=Ru$
\item[2.]if
$\nu_y(\bar{u})\leq\varepsilon$, then $u^*=(x^*;\l^*):=(\bar{x};\bar{\l})$ else;
\item[3.] find  $\bar{\tau}=\max\{0\leq\tau\leq 1 :x(\tau)=y+t(\bar{x}-y)\in\Omega_{\gamma},\;(\nabla f(x(\tau)),\bar{x}-y)\leq 0\}$ and $x(\bar{\tau})$;
\item[4.] if
\begin{equation}\label{Eq3_27.11.16}
f(y)-f(x(\bar{\tau}))\geq \Delta,
\end{equation}
 then update the center $\bar{y}:=0.5(y+x(\bar{x}))$,  set $x:=\bar{x}$, $y:=\bar{y}$ and go to 1;\\
 else  set $x:=\bar{x};$ $\l:=\bar{\l}$ and go to 1.
 \end{enumerate}
It follows from 3. and 4.  that the sequence of centers is monotone decreasing in value, therefore from some point on the inequality \eqref{Eq3_27.11.16} can't be satisfied, so the "center" is fixed. 

Hence, from this point on the primal-dual sequence is generated only by the relaxation operator (\ref{Eq4_27.11.16}) and converge to the primal-dual solution with $Q$-linear rate due to Theorem \ref{T_27.11.16}

%--------------------------------------
\section{Concluding Remarks}
It follows from 1.-4. that the efficiency  of the EPM heavily depends on the efficiency  of the unconstrained minimization algorithm used in operator $R$.

The absence of singularity of EDF  at the solution combined with stability of its Hessian's condition number improves substantially the  efficiency of the operator $R$.In particular, it allows to reduce the number of unconstrained minimization steps per Lagrange multipliers update.

On the other hand, under  fixed  Lagrange multipliers EDF possesses self-concordance properties    for a wide classes of constrained optimization problems.

It provides an opportunity to combine the nice feature of the IPM at the beginning of the computational process with excellent EDF properties at the final phase.

The NR approach  produced very strong  numerical results for  wide classes of large scale nonlinear optimization problems (see, for example, \cite{Ben-Tal97}, \cite{Breirfeld96}, \cite{Griva06}, \cite{Nash94}).

 In particular,  one of the most reliable NLP solver PENNON is based on NR theory (see \cite{Kocvara03}, \cite{Kocvara05}).

It leads us to believe that the extra tool, which EDF possesses, can contribute to the numerical efficiency  mainly because updating the center does not require much computational effort, but can substantial reduce the  objective function value.

It means that updating the "center" will allow to reach the "hot start" faster (see \cite{Griva06}, \cite{polyak92}, \cite{polyak16}).

\end{document}